# New $M$-estimators in semi-parametric regression with errors in variables

Cristina Butucea[a,b] and Marie-Luce Taupin[c,d]

[a]Université Paris X, Modal'X, 200, ave. de la République, 92001 Nanterre Cedex, France
[b]Université Paris VI, PMA, 175, rue de Chevaleret, 75013 Paris, France. E-mail: butucea@ccr.jussieu.fr
[c]Université Paris Descartes, IUT de Paris, 143 ave. de Versailles, 75016 Paris, France
[d]Université Paris-Sud, Bât. 425, Département de Mathématiques, 91405 Orsay Cedex, France.
E-mail: marie-luce.taupin@math.u-psud.fr

**Abstract.** In the regression model with errors in variables, we observe $n$ i.i.d. copies of $(Y, Z)$ satisfying $Y = f_{\theta^0}(X) + \xi$ and $Z = X + \varepsilon$ involving independent and unobserved random variables $X, \xi, \varepsilon$ plus a regression function $f_{\theta^0}$, known up to a finite dimensional $\theta^0$. The common densities of the $X_i$'s and of the $\xi_i$'s are unknown, whereas the distribution of $\varepsilon$ is completely known. We aim at estimating the parameter $\theta^0$ by using the observations $(Y_1, Z_1), \ldots, (Y_n, Z_n)$. We propose an estimation procedure based on the least square criterion $\tilde{S}_{\theta^0, g}(\theta) = \mathbb{E}_{\theta^0, g}[((Y - f_\theta(X))^2 w(X)]$ where $w$ is a weight function to be chosen. We propose an estimator and derive an upper bound for its risk that depends on the smoothness of the errors density $p_\varepsilon$ and on the smoothness properties of $w(x) f_\theta(x)$. Furthermore, we give sufficient conditions that ensure that the parametric rate of convergence is achieved. We provide practical recipes for the choice of $w$ in the case of nonlinear regression functions which are smooth on pieces allowing to gain in the order of the rate of convergence, up to the parametric rate in some cases. We also consider extensions of the estimation procedure, in particular, when a choice of $w_\theta$ depending on $\theta$ would be more appropriate.

**Résumé.** Dans le modèle de régression avec erreurs sur les variables, nous observons $n$ v.a. i.i.d. de même loi que $(Y, Z)$ satisfaisant aux relations $Y = f_{\theta^0}(X) + \xi$ et $Z = X + \varepsilon$, où les v.a. $X, \xi, \varepsilon$ sont indépendantes, pas observées, et la fonction de régression $f_{\theta^0}$ est connue à un paramètre de dimension finie $\theta^0$ près. Les densités de $X$ et de $\xi$ sont inconnues tandis que la loi de $\varepsilon$ est entièrement connue. Nous estimons le paramètre $\theta^0$ à partir des observations $(Y_1, Z_1), \ldots, (Y_n, Z_n)$. Nous proposons une procédure d'estimation basée sur le critère des moindres carrés $\tilde{S}_{\theta^0, g}(\theta) = \mathbb{E}_{\theta^0, g}[((Y - f_\theta(X))^2 w(X)]$, où $w$ est une fonction de poids à choisir. Nous définissons l'estimateur et calculons la borne supérieure du risque de cet estimateur, qui dépend de la régularité de la densité des erreurs $p_\varepsilon$ et de la régularité en $x$ de $w(x) f_\theta(x)$. De plus, nous établissons des conditions suffisantes pour que les estimateurs atteignent la vitesse paramétrique. Nous décrivons des méthodes pratiques pour le choix de $x$ dans le cas des fonctions de régression non-linéaires qui sont régulières par morceaux permettant de gagner des ordres de vitesse allant jusqu'à la vitesse paramétrique dans certains cas. Nous considérons également des extensions de cette procédure d'estimation, en particulier au cas où un choix de $w_\theta$ dépendant de $\theta$ serait plus appropié.

## 1. Introduction

We consider the regression model with errors in variables where one observes $n$ independent and identically distributed (i.i.d.) copies of $(Y, Z)$ satisfying

$$\begin{cases} Y = f_{\theta^0}(X) + \xi, \\ Z = X + \varepsilon, \end{cases}$$

involving independent and unobserved random variables $X, \xi, \varepsilon$, plus a regression function $f_{\theta^0}$ known up to a finite dimensional parameter $\theta^0$, belonging to the interior of a compact set of $\Theta \subset \mathbb{R}^d$. The common densities of the $X_i$'s and of the $\xi_i$'s are unknown, with $\mathbb{E}(\xi) = 0$, whereas the density of the errors $\varepsilon_i$'s is completely known. In this context, we aim at estimating the finite dimensional parameter $\theta^0$ in the presence of a functional nuisance parameter $g$, the density of the $X$.

*Previous known results*

This model has been widely studied with first results written in the 1950's (see, for instance, [28] or [20]). Most of the results deal with linear models where $\sqrt{n}$-consistency, asymptotic normality, and efficiency have been studied. One can cite among the others Bickel and Ritov [3], Bickel et al. [2], Cheng and van Ness [7], van der Vaart [30, 31, 32], Murphy and van der Vaart [26]. The nonlinear models have been more recently considered starting with the case of repeated measurement data as in [12, 33, 34], under additional assumptions as in [13, 17, 21, 22, 25] or by simulation (see [5, 18, 19, 24]). To our knowledge, the first consistent estimator in nonlinear regression models with errors in variables, under nonparametric assumptions on the design density $g$ has been proposed by Taupin [29] when the errors $\varepsilon$ are Gaussian and by Comte and Taupin [8] in the context of auto-regressive models with errors in variables for various types of errors density $p_\varepsilon$. In those papers, the estimation procedure is based on the estimation of the modified least square criterion $\mathbb{E}[(Y - \mathbb{E}(f_\theta(X)|Z))^2 W(Z)]$ where the conditional expectation is estimated by using the observations $Z_1, \ldots, Z_n$ and where $W$ is a compactly supported weight function. It is also stated that the rate of convergence, which does not have an explicit form, depends on the smoothness of the regression function as well as the smoothness of $p_\varepsilon$ through the increase of the ratio $(f_\theta p_\varepsilon(z - \cdot))^*(t)/p_\varepsilon^*(t)$ as $t$ goes to infinity. The parametric rate of convergence is achieved in some specific examples, such as polynomial or exponential regression functions. The main drawback of this estimator is that besides its complexity, its rate of convergence does not have an explicit form.

More recently, Hong and Tamer [16] propose a consistent estimator in the specific case where $p_\varepsilon^*$ is of the form $p_\varepsilon^*(t) = 1/(1 + \sigma^2 t^2)$. Their estimation procedure strongly depends on this particular form of $p_\varepsilon^*$ through the ratio $1/p_\varepsilon^*$ always appearing in errors in variables techniques. The extension of their method for other errors density $p_\varepsilon$ seems thus not obvious.

*Our results*

We propose here a new estimation procedure based on the least square criterion that is more general, explicit, natural, and tractable than the one proposed by Taupin [29]. Moreover, it often provides better results and it allows to provide sufficient conditions to achieve the parametric rate.

The procedure is based on the estimation, using the observations $(Y_i, Z_i)$ for $i = 1, \ldots, n$, of the least squares criterion $\tilde{S}_{\theta^0,g} := \mathbb{E}_{\theta^0,g}[(Y - f_\theta(X))^2 w(X)]$ where $w$ is a positive weight function to be chosen.

In this context, we naturally define in Section 2 the estimator

$$\tilde{\theta}_1 = \arg\min_{\theta \in \Theta} \tilde{S}_{n,1}(\theta) \quad \text{with } \tilde{S}_{n,1}(\theta) = \frac{1}{n} \sum_{i=1}^n \int [(Y_i - f_\theta(x))^2 w(x)] C_n K_n(C_n(x - Z_i)) \, dx,$$

where $K_n$ is a kernel such that its Fourier transform satisfies $K_n^*(t) = K^*(t)/p_\varepsilon^*(tC_n)$, for $K^*$ compactly supported. In the literature, such a kernel $K_n$ is commonly known as the deconvolution kernel. In the sequel, $C_n$ is a sequence which tends to $\infty$ and $p^*$ denotes the Fourier transform, defined by $p^*(t) = \int e^{itx} p(x) \, dx$ of an arbitrary square integrable function $p$.

We show in Section 3 that under classical identifiability, moment and smoothness assumptions, $\tilde{\theta}_1$ is a consistent estimator of $\theta^0$ with a rate of convergence depending on two factors, the smoothness of $(wf_\theta)$, and the smoothness of the errors density $p_\varepsilon$. This partly comes from the fact that this estimation procedure is based on the estimation of the linear integral functional $\mathbb{E}_{\theta^0,g}(Yw(X)f_\theta(X))$ and $\mathbb{E}_{\theta^0,g}(w(X)f_\theta^2(X))$ using the observations $(Y_1, Z_1), \ldots, (Y_n, Z_n)$, that is, by recovering information on $(Y, X)$ using the observation $(Y, Z)$. More precisely, it depends on the smoothness of $w(x)\partial f_\theta(x)/\partial \theta$ and $w(x)\partial(f_\theta^2(x))/\partial \theta$ and the smoothness of the errors density $p_\varepsilon(x)$, as a function of $x$ through the behavior of $(w\partial f_\theta/\partial \theta)^*(t)/p_\varepsilon^*(t)$, $(w\partial(f_\theta^2)/\partial \theta)^*(t)/p_\varepsilon^*(t)$ as function of $t \to \infty$. From this construction, we derive sufficient conditions that ensure that $\tilde{\theta}_1$ achieves the parametric rate of convergence.

The rate of convergence of the proposed estimator is thoroughly studied for various smoothness properties of $wf_\theta$, $wf_\theta^2$ as well as their derivatives in $\theta$ as functions of $x$ and for various errors' densities $p_\varepsilon$. It appears that as for the nonparametric estimation of the regression function in errors in variables models, the smoother is the noise density $p_\varepsilon$ the slower becomes the rate of convergence of estimators. Nevertheless, in examples we can considerably improve the smoothness of the functions $f_\theta$ and $f_\theta^2$ by multiplying them by a properly chosen weight function. This automatically improves the rates of convergence of our estimator.

The conditions ensuring that the $\sqrt{n}$-consistency is achieved are also deeply studied. The main idea is that the actual shape of the regression function matters less than its smoothness (compared to the noise smoothness). Those conditions are illustrated through various examples of regression functions in Section 4. The point is that these conditions allow to achieve $\sqrt{n}$-consistency in setups that were not known before.

This estimation procedure is used and extended to construct various other estimators.

Firstly, under conditions ensuring that $\mathbb{E}_{\theta^0,g}(Yw(X)f_\theta(X))$ and $\mathbb{E}_{\theta^0,g}(w(X)f_\theta^2(X))$ are estimated at the parametric rate $\sqrt{n}$, we propose in Section 5 a second estimator $\tilde{\theta}_2$ of $\theta^0$. In that case, the parametric rate of convergence for the estimation of $\theta^0$ can be achieved and the asymptotic normality of this estimator is stated. The first estimator, $\tilde{\theta}_1$, based on a deconvolution kernel is more general and applicable in all setups, but for specific regression functions, the use of a deconvolution kernel is not required. In those cases, the second estimator is more simple. Nevertheless, the conditions ensuring that $\mathbb{E}_{\theta^0,g}(Yw(X)f_\theta(X))$ and $\mathbb{E}_{\theta^0,g}(w(X)f_\theta^2(X))$ are estimated at the parametric rate $\sqrt{n}$ are not always fulfilled.

Secondly, when the variance $\sigma_{\xi,2}^2 = \text{Var}(\xi)$ is known, we propose in Section 6 another estimation criterion based on $w_\theta$ depending on $\theta$. It allows to improve the rate of convergence of the estimators, by smoothing in a better manner $w_\theta f_\theta$. For a large class of regression functions, $w$ not depending on $\theta$ works well, for instance for polynomial, exponential, cosines regression functions as well as for regression functions $f_\theta$ of the form $f_\theta(x) = \varphi(\theta)f(x)$. But sometimes, the smoothing properties of $w$ will be improved by taking $w$ depending on $\theta$, e.g., in the case where the regression function has to be smoothed at some point related to $\theta$.

In this context, we extend our procedure to the estimation of the least squares criterion

$$S_{\theta^0,g} := \mathbb{E}_{\theta^0,g}[(Y - f_\theta(X))^2 w_\theta(X)] - \sigma_{\xi,2}^2 \mathbb{E}_{\theta^0,g}[w_\theta(X)],$$

using the observations $(Y_1, Z_1), \ldots, (Y_n, Z_n)$, where $w_\theta$ is a positive weight function, depending on $\theta$, to be chosen and where $\sigma_{\xi,2}^2 = \text{Var}(\xi)$ is now supposed to be known. Using this criterion and analogously to the construction of $\tilde{\theta}_1$, we propose to estimate $\theta^0$ by $\widehat{\theta}_1$ defined by

$$\widehat{\theta}_1 = \arg\min_{\theta \in \Theta} S_{n,1}(\theta) \quad \text{with } S_{n,1}(\theta) = \frac{1}{n} \sum_{i=1}^n \int [(Y_i - f_\theta(x))^2 - \sigma_{\xi,2}^2] w_\theta(x) C_n K_n(C_n(x - Z_i)) \, dx.$$

Under classical identifiability, moment and smoothness assumptions, the estimator $\widehat{\theta}_1$ is a consistent estimator of $\theta^0$ with a rate of convergence depending on the smoothness of $\partial(w_\theta(x)f_\theta(x))/\partial \theta$, $\partial(w_\theta(x)f_\theta^2(x))/\partial \theta$ and on the smoothness of the errors density $p_\varepsilon(x)$, as functions of $x$. From this construction, we derive sufficient conditions that ensure that $\widehat{\theta}_1$ achieves the parametric rate of convergence.

As for the first extension, we propose another estimator $\hat{\theta}_2$ of $\theta^0$ based on $S_{\theta^0,g}$ achieving $\sqrt{n}$-consistency under constructive conditions ensuring that $\mathbb{E}_{\theta^0,g}(Y^2 w_\theta(X))$, $\mathbb{E}_{\theta^0,g}(Y w_\theta(X) f_\theta(X))$ and $\mathbb{E}_{\theta^0,g}(w_\theta(X) f_\theta^2(X))$ can be estimated at the parametric rate $\sqrt{n}$.

Let us have a look to the properties of our estimator in the context studied by Hong and Tamer [16]. Our estimation procedure, more general, allows to recover $\sqrt{n}$-consistency, in the case of noise densities satisfying $p_\varepsilon^*(t) = c|t|^{-2}(1 + \mathrm{o}(1))$ as $|t| \to \infty$ with regression functions $f_\theta$ having derivatives in $\theta$ up to order 3, twice continuously differentiable functions of $x$.

The drawback of smoothing by multiplication is that for particular regression functions, we can obtain infinitely differentiable functions but not analytic ones. In such a case, and if the noise has an analytic density, the parametric rate of convergence can not be attained by our method. Nevertheless, the smoothing technique significantly improves the rate by using a clever choice of weight function $w$.

The main question that remains open is: Is it possible to construct a $\sqrt{n}$-consistent estimator of $\theta^0$ for all regression functions and without any restriction on the smoothness of the errors density?

The paper is organized as follows. Section 2 presents the estimation procedure and Section 3 gives the asymptotic properties of this estimator. Those asymptotic properties and practical recipes are illustrated through examples in Section 4. Section 5 presents another estimator and Section 6 gives extensions of the two former estimators, illustrated by examples in Section 7.

The proofs can be found in Section 8 with technical lemmas presented in the Appendix.

## 2. The estimation procedure

*Notations*

We denote by $x_-$ the negative part of $x$, $\| \varphi \|_2^2 = \int \varphi^2(x)\,\mathrm{d}x$, $\| \varphi \|_\infty = \sup_{x \in \mathbb{R}} |\varphi(x)|$. In the same way, $p \star q(z) = \int p(z - x)q(x)\,\mathrm{d}x$ denotes the convolution of two square integrable functions $p$ and $q$. The variance $\mathrm{Var}(\xi) = \mathbb{E}(\xi^2)$ is denoted by $\sigma_{\xi,2}^2$ and $\mathbb{E}(\xi^4) = \sigma_{\xi,4}^4$. For $\theta \in \mathbb{R}^d$, $\| \theta \|_{\ell^2}^2 = \sum_{k=1}^d \theta_k^2$ and $\theta^\top$ is the transpose matrix of $\theta$.

From now, $\mathbb{P}$, $\mathbb{E}$, and Var denote the probability $\mathbb{P}_{\theta^0,g}$, the expected value $\mathbb{E}_{\theta^0,g}$ and, respectively, the variance $\mathrm{Var}_{\theta^0,g}$, when the underlying and unknown true parameters are $\theta^0$ and $g$.

The starting point of the estimation procedure is to construct an estimator based on the observations $(Y_i, Z_i)$ for $i = 1, \ldots, n$, of the least square criterion

$$\tilde{S}_{\theta^0,g}(\theta) = \mathbb{E}[(Y - f_\theta(X))^2 w(X)], \tag{2.1}$$

where $w$ is a positive weight function to be suitably chosen.

This estimation procedure requires the following assumptions.

*Smoothness and moment assumptions*

(A$_1$) For any $\theta$ in $\Theta$, the function $\theta \mapsto f_\theta$ admits continuous derivatives with respect to $\theta$ up to the order 3.
(A$_2$) The quantities $\mathbb{E}[w^2(X)(Y - f_\theta(X))^4]$ and their derivatives with respect to $\theta$ up to order 2 are finite.

We denote by $\mathcal{G}$ the set of densities $g$ such that (A$_2$) holds. We subsequently assume that there exist two constants $C(f_{\theta^0}^2)$ and $C(f_{\theta^0})$, depending only on $\theta^0$ through the functions $f_{\theta^0}^2$ and $f_{\theta^0}$, respectively, such that

(A$_3$) $\sup_{g \in \mathcal{G}} \| f_{\theta^0}^2 g \|_2^2 \leq C(f_{\theta^0}^2)$, $\sup_{g \in \mathcal{G}} \| f_{\theta^0} g \|_2^2 \leq C(f_{\theta^0})$.
(A$_4$) $\sup_{\theta \in \Theta} |wf_\theta|$, $|w|$ and $\sup_{\theta \in \Theta} |wf_\theta^2|$ belong to $\mathbb{L}_1(\mathbb{R})$.

*Identifiability assumptions*

(II$_1$) The quantity $\tilde{S}_{\theta^0,g}(\theta) = \sigma_{\xi,2}^2 \mathbb{E}(w(X)) + \mathbb{E}[(f_{\theta^0}(X) - f_\theta(X))^2 w(X)]$ admits one unique minimum at $\theta = \theta^0$.

(II$_2$) For all $\theta \in \Theta$ the matrix $\tilde{S}^{(2)}_{\theta^0,g}(\theta) = \partial^2 S_{\theta^0,g}(\theta)/\partial\theta^2$ exists and the matrix

$$\tilde{S}^{(2)}_{\theta^0,g}(\theta^0) = \mathbb{E}\left[2w(X)\left(\frac{\partial f_\theta(X)}{\partial \theta}\bigg|_{\theta=\theta^0}\right)\left(\frac{\partial f_\theta(X)}{\partial \theta}\bigg|_{\theta=\theta^0}\right)^\top\right] \qquad \text{is positive definite.}$$

### 2.1. Construction of the estimator

We denote by $f_{Y,X}$ and $f_{Y,Z}$ the joint densities of $(Y,X)$ and respectively of $(Y,Z)$ satisfying in this model,

$$f_{Y,X}(y,x) = g(x)p_\xi(y - f_{\theta^0}(x)) \quad \text{and} \quad f_{Y,Z}(y,z) = f_{Y,X}(y,\cdot) \star p_\varepsilon(z). \tag{2.2}$$

Write $\tilde{S}_{\theta^0,g}(\theta)$ as

$$\tilde{S}_{\theta^0,g}(\theta) = \mathbb{E}[(Y - f_\theta(X))^2 w(X)] = \int (y - f_\theta(x))^2 w(x) f_{Y,X}(y,x) \, dy \, dx.$$

We consider $p_\varepsilon$ that satisfies the following assumption.

(N$_1$) The density $p_\varepsilon$ belongs to $\mathbb{L}_2(\mathbb{R})$ and for all $x \in \mathbb{R}, p_\varepsilon^*(x) \neq 0$.

The assumption (N$_1$) is quite usual in density deconvolution.

According to (2.2) and under (N$_1$), we naturally propose to estimate $\tilde{S}_{\theta^0,g}(\theta)$ by

$$\tilde{S}_{n,1}(\theta) = \frac{1}{n}\sum_{i=1}^n \int (Y_i - f_\theta(x))^2 w(x) K_{n,C_n}(x - Z_i) \, dx = \frac{1}{n}\sum_{i=1}^n ((Y_i - f_\theta)^2 w) \star K_{n,C_n}(Z_i), \tag{2.3}$$

where $K_{n,C_n}(\cdot) = C_n K_n(C_n \cdot)$ is a deconvolution kernel defined via its Fourier transform, such that $\int K_n(x) \, dx = 1$ and

$$K^*_{n,C_n}(t) = \frac{K^*_{C_n}(t)}{p_\varepsilon^*(t)} = \frac{K^*(t/C_n)}{p_\varepsilon^*(t)}, \tag{2.4}$$

with $K^*$ compactly supported satisfying $|1 - K^*(t)| \leq \mathbb{1}_{|t|\geq 1}$.

Using this empirical criterion, we propose to estimate $\theta^0$ by

$$\tilde{\theta}_1 = \arg\min_{\theta \in \Theta} \tilde{S}_{n,1}(\theta). \tag{2.5}$$

## 3. Asymptotic properties

### 3.1. General results for the risk of $\tilde{\theta}_1$

We say that a function $\psi \in \mathbb{L}_1(\mathbb{R})$ satisfies (3.6) if for a sequence $C_n$ and under previous notation we have

$$\min_{q=1,2} \|\psi^*(K^*_{C_n} - 1)\|_q^2 + n^{-1} \min_{q=1,2} \left\|\frac{\psi^* K^*_{C_n}}{p_\varepsilon^*}\right\|_q^2 = \mathrm{o}(1). \tag{3.6}$$

**Theorem 3.1.** *Let $\tilde{\theta}_1 = \tilde{\theta}_1(C_n)$ be defined by (2.5) under the assumptions (II$_1$), (II$_2$), (N$_1$), (A$_1$)–(A$_3$) and (A$_4$).*

(1) *For all of the sequences $C_n$ such that $w$, $wf_\theta$ and $wf_\theta^2$ and their first derivatives with respect to $\theta$ satisfy (3.6), $\mathbb{E}(\|\tilde{\theta}_1 - \theta^0\|_{\ell^2}^2) = \mathrm{o}(1)$, as $n \to \infty$ and $\tilde{\theta}_1 = \tilde{\theta}_1(C_n)$ is a consistent estimator of $\theta^0$.*

(2) *Assume, moreover, that for all $\theta \in \Theta$, $w$, $f_\theta w$ and $f_\theta^2 w$ and their derivatives up to order 3 with respect to $\theta$ satisfy* (3.6). *Then* $\mathbb{E}(\|\tilde{\theta}_1 - \theta^0\|_{\ell^2}^2) = \mathrm{O}(\tilde{\varphi}_n^2)$ *with* $\tilde{\varphi}_n = \|(\tilde{\varphi}_{n,j})\|_{\ell^2}$, $\tilde{\varphi}_{n,j}^2 = \tilde{B}_{n,j}^2(\theta^0) + \tilde{V}_{n,j}(\theta^0)/n$, $j = 1, \ldots, d$, *where*

$$\tilde{B}_{n,j}(\theta) = \min\{\tilde{B}_{n,j}^{[1]}(\theta), B_{n,j}^{[2]}(\theta)\} \quad \text{and} \quad \tilde{V}_{n,j}(\theta) = \min\{\tilde{V}_{n,j}^{[1]}(\theta), \tilde{V}_{n,j}^{[2]}(\theta)\}$$

*and for* $q = 1, 2$

$$\tilde{B}_{n,j}^{[q]}(\theta) = \left\|\left(\frac{\partial(wf_\theta)}{\partial \theta_j}\right)^* (K_{C_n}^* - 1)\right\|_q^2 + \left\|\left(\frac{\partial(wf_\theta^2)}{\partial \theta_j}\right)^* (K_{C_n}^* - 1)\right\|_q^2$$

*and*

$$\tilde{V}_{n,j}^{[q]}(\theta) = \left\|\left(\frac{\partial(wf_\theta)}{\partial \theta_j}\right)^* \frac{K_{C_n}^*}{p_\varepsilon^*}\right\|_q^2 + \left\|\left(\frac{\partial(wf_\theta^2)}{\partial \theta_j}\right)^* \frac{K_{C_n}^*}{p_\varepsilon^*}\right\|_q^2.$$

**Remark 3.1.** *It is noteworthy that for any integrable $\psi$, one can always find a sequence $C_n$ such that (3.6) holds.*

**Remark 3.2.** *The terms $\tilde{B}_{n,j}^2$ and $\tilde{V}_{n,j}/n$ are, respectively, the squared bias and variance terms with as in density deconvolution, bigger variance for smoother error density $p_\varepsilon$ and smaller bias for smoother $(wf_\theta)$.*

As we can see, the rate of convergence for estimating $\theta^0$ is given by both the smoothness of the errors density $p_\varepsilon$ and the smoothness of $f_\theta w$; more precisely by the smoothness of $\partial(f_\theta w)/\partial \theta$ and $\partial(f_\theta^2 w)/\partial \theta$, as functions of $x$. Those smoothness properties are described in both cases, by the asymptotic behavior of the Fourier transforms. The slower rates are obtained for the smoother errors density $p_\varepsilon$, for instance for Gaussian $\varepsilon$'s. Nevertheless, those rates are improved with a proper choice of $w$ such that $(wf_\theta)$, $(wf_\theta^2)$ and their derivatives with respect to $\theta$ are smooth enough. In this context, the parametric rate of convergence is achieved as soon as the derivatives with respect to $\theta$, of $(wf_\theta)$ and $(wf_\theta^2)$ as functions of $x$ are smoother than the errors density $p_\varepsilon$.

### 3.2. Consequence: sufficient conditions to obtain the parametric rate of convergence

We say that the conditions $(C_1)$–$(C_3)$ hold if there exists a weight function $w$ such that for all $\theta \in \Theta$,

($C_1$) the functions $(wf_\theta)$ and $(wf_\theta^2)$ belong to $\mathbb{L}_1(\mathbb{R})$, and the functions $\sup_\theta w^*/p_\varepsilon^*$, $\sup_\theta (f_\theta w)^*/p_\varepsilon^*$, $\sup_\theta (f_\theta^2 w)^*/p_\varepsilon^*$ belong to $\mathbb{L}_1(\mathbb{R}) \cap \mathbb{L}_2(\mathbb{R})$;

($C_2$) the functions $\sup_{\theta \in \Theta}(\frac{\partial(f_\theta w)}{\partial \theta})^*/p_\varepsilon^*$ and $\sup_{\theta \in \Theta}(\frac{\partial(f_\theta^2 w)}{\partial \theta})^*/p_\varepsilon^*$ belong to $\mathbb{L}_1(\mathbb{R}) \cap \mathbb{L}_2(\mathbb{R})$;

($C_3$) the functions $(\frac{\partial^2(f_\theta w)}{\partial \theta^2})^*/p_\varepsilon^*$ and $(\frac{\partial^2(f_\theta^2 w)}{\partial \theta^2})^*/p_\varepsilon^*$ belong to $\mathbb{L}_1(\mathbb{R}) \cap \mathbb{L}_2(\mathbb{R})$.

**Theorem 3.2.** *Consider model (1.1) under the assumptions $(A_1)$–$(A_4)$, $(II_1)$, $(II_2)$, $(N_1)$ and the conditions $(C_1)$–$(C_3)$. Consider $C_n$ such that for all $\theta \in \Theta$, $f_\theta w$, $f_\theta^2 w$ and their derivatives up to order 3 satisfy (3.6). Then $\tilde{\theta}_1$ defined by (2.5) is a $\sqrt{n}$-consistent estimator of $\theta^0$ which satisfies moreover that $\sqrt{n}(\tilde{\theta}_1 - \theta^0) \xrightarrow[n \to \infty]{\mathcal{L}} \mathcal{N}(0, \tilde{\Sigma}_1)$, with $\tilde{\Sigma}_1$ that equals*

$$\left(\mathbb{E}\left[2w(X)\left(\frac{\partial f_\theta(X)}{\partial \theta}\right)\left(\frac{\partial f_\theta(X)}{\partial \theta}\right)^\top\right]\bigg|_{\theta=\theta^0}\right)^{-1} \tilde{\Sigma}_{0,1} \left(\mathbb{E}\left[2w(X)\left(\frac{\partial f_\theta(X)}{\partial \theta}\right)\left(\frac{\partial f_\theta(X)}{\partial \theta}\right)^\top\right]\bigg|_{\theta=\theta^0}\right)^{-1},$$

*where $\tilde{\Sigma}_{0,1}$ is given by*

$$(2\pi)^{-2}\mathbb{E}\left\{\left[\int \left(\frac{\partial(f_\theta^2 w - 2Y f_\theta w)}{\partial \theta}\bigg|_{\theta=\theta^0}\right)^*(u) \frac{e^{-iuZ}}{p_\varepsilon^*(u)} du\right]\left[\int \left(\frac{\partial(f_\theta^2 w - 2Y f_\theta w)}{\partial \theta}\bigg|_{\theta=\theta^0}\right)^*(u) \frac{e^{-iuZ}}{p_\varepsilon^*(u)} du\right]^\top\right\}.$$

### 3.3. Resulting rates for general smoothness classes

We now precise the asymptotic properties of $\tilde{\theta}_1$ when the decrease of these Fourier transforms are quantified through the following assumptions.

(N$_2$) There exist positive constants $\underline{C}(p_\varepsilon), \overline{C}(p_\varepsilon), \beta, \rho, \alpha,$ and $u_0$ such that

$$\underline{C}(p_\varepsilon) \leq |p_\varepsilon^*(u)||u|^\alpha \exp(\beta|u|^\rho) \leq \overline{C}(p_\varepsilon) \quad \text{for all } |u| \geq u_0.$$

If $\rho > 0$ in (N$_2$), then the noise is called exponential noise or super smooth noise. And if $\rho = 0$ in (N$_2$), then by convention $\beta = 0$, $\alpha > 1$ and the noise is called polynomial noise or ordinary smooth.

The smoothness properties of functions involving the regression function are also given by the asymptotic behavior of the Fourier transforms described as follows.

(R$_1$) A function $f$ satisfies (R$_1$) if $f$ belongs to $\mathbb{L}_1(\mathbb{R}) \cap \mathbb{L}_2(\mathbb{R})$ and if there exist $a, b, r,$ and $u_0 \geq 0$ such that
$$\underline{L}(f) \leq |f^*(u)||u|^a \exp(b|u|^r) \leq \overline{L}(f) < \infty \text{ for all } |u| \geq u_0.$$

If $r = 0$ in (R$_1$) then by convention $b = 0$ and the function $f$ is called ordinary smooth. If $r > 0$, the function $f$ is called super smooth.

**Corollary 3.1.** *Under the assumptions of Theorem* 3.1, *assume that $p_\varepsilon$ satisfies* (N$_2$) *and that for all $\theta \in \Theta$, $(f_\theta w)$, $(f_\theta^2 w)$ and their derivatives with respect to $\theta_j$, $j = 1, \ldots, d$ up to order 3, satisfy* (R$_1$).

(1) *For all the sequences $C_n$ such that*

$$C_n^{(2\alpha - 2a + 1 - \rho)} \exp\{-2bC_n^r + 2\beta C_n^\rho\}/n = o(1) \quad \text{as } n \to +\infty, \tag{3.7}$$

$\mathbb{E}(\|\tilde{\theta}_1 - \theta^0\|_{\ell^2}^2) = O(\tilde{\varphi}_n^2) = o(1)$ *as* $n \to \infty$ *and* $\tilde{\theta}_1 = \tilde{\theta}_1(C_n)$ *is a consistent estimator of $\theta^0$.*
(2) *Moreover, $\tilde{\varphi}_n^2$ is given by Table* 1, *according to values of parameters $a, b, r, \alpha, \beta$ and $\rho$.*

## 4. Examples and methodological advice (1)

In this section, we propose a deep study of the asymptotic properties of the estimator through various examples of regression functions. We show that for many regression functions the practitioner may encounter

Table 1
Rates of convergence $\tilde{\varphi}_n^2$ of $\tilde{\theta}_1$

| | | $p_\varepsilon$ | |
|---|---|---|---|
| | | $\rho = 0$ in (N$_2$) ordinary smooth | $\rho > 0$ in (N$_2$) super smooth |
| $wf_{\theta^0}$ | (R$_1$) $b = r = 0$ Sobolev | $a < \alpha + 1/2$ $\quad n^{-(2a-1)/(2\alpha)}$ $a \geq \alpha + 1/2$ $\quad n^{-1}$ | $(\log n)^{-(2a-1)/\rho}$ |
| | (R$_1$) $r > 0$ | $n^{-1}$ | $r < \rho$ $\quad (\log n)^{A(a,r,\rho)} \exp\{-2b(\frac{\log n}{2\beta})^{r/\rho}\}$ |
| | $\mathcal{C}^\infty$ | | $r = \rho$ $\quad b < \beta$ $\quad \frac{(\log n)^{A(a,r,\rho) + 2\alpha b/(\beta r)}}{n^{b/\beta}}$ |
| | | | $\quad\quad\quad\quad b = \beta, a < \alpha + 1/2$ $\quad \frac{(\log n)^{(2\alpha - 2a + 1)/r}}{n}$ |
| | | | $\quad\quad\quad\quad b = \beta, a \geq \alpha + 1/2$ $\quad n^{-1}$ |
| | | | $\quad\quad\quad\quad b > \beta$ $\quad n^{-1}$ |
| | | | $r > \rho$ $\quad n^{-1}$ |

Where $A(a, r, \rho) = (-2a + 1 - r + (1 - r)_-)/\rho$.

there are a few simple smoothing weight functions to choose so that the rates improve significantly, up to parametric rate in many cases. This new procedure allows to achieve the parametric rate of convergence in lot of examples and especially in examples where the previously known estimator proposed in [29] does not. In all of these examples, the noise distribution is arbitrary, as far as it satisfies $(N_1)$ and $(N_2)$ with $\rho \leq 2$. The two first examples simply show that this new and more general estimation procedure allows to recover previous known results in simple cases. The other examples provide new results that underlie the improvement due to this method.

*Example 1 (Polynomial regression function).* Let $f_\theta$ be of the form $f_\theta(x) = \sum_{k=1}^p \theta_k x^k$ and $p_\varepsilon$ satisfying $(N_2)$ with $\rho \leq 2$. Assume that $\mathbb{E}(Y^2) < \infty$ and that $\mathbb{E}(Z^{2p}) < \infty$. Let $K$ be such $K^*(t) = \mathbb{1}_{|t| \leq 1}$ and let $w(x) = \exp\{-x^2/(4\beta)\}$. Then conditions $(C_1)$–$(C_3)$ are satisfied. Consequently, the estimator $\tilde{\theta}_1$ is a $\sqrt{n}$-consistent and asymptotically Gaussian estimator of $\theta^0$.

*Remark 4.1.* In this example, one can also choose $w \equiv 1$, provided that the kernel $K$ has finite absolute moments of order $p$ and satisfies $\int u^r K(u)\,\mathrm{d}u = 0$, for $r = 1, \ldots, p$. With these choices of $w$ and $K$, $\tilde{\theta}_1$ remains a $\sqrt{n}$-consistent and asymptotically Gaussian estimator of $\theta^0$.

The $\sqrt{n}$-consistency as well as the asymptotic normality was already achieved with different estimators, in the linear case (see, e.g., [2, 3, 26, 31, 32]). In polynomial case, other $\sqrt{n}$-consistent estimators already exist, without proving the asymptotic normality (see, e.g. [29] and [8]) or in the polynomial functional errors in variables model, with fixed and not random $X_i$'s (see, e.g., [14, 15] or [6]). It is noteworthy that our new estimation procedure, quite more simple and natural than the one proposed in [29], also provides the $\sqrt{n}$-consistency in this simple case.

*Example 2 (Exponential regression function).* Let $f_\theta$ be of the form $f_\theta(x) = \exp(\theta x)$ and $p_\varepsilon$ satisfying $(N_2)$ with $\rho \leq 2$. Assume that $\mathbb{E}(Y^2) < \infty$ and that $\mathbb{E}[\exp(2\theta^0 Z)] < \infty$. Let $K$ be such $K^*(t) = \mathbb{1}_{|t| \leq 1}$ and let $w(x) = \exp\{-x^2/(4\beta)\}$. Then the conditions $(C_1)$–$(C_3)$ are satisfied and the estimator $\tilde{\theta}_1$ is a $\sqrt{n}$-consistent and asymptotically Gaussian estimator of $\theta^0$. Once again, this new estimation procedure allows to achieve the $\sqrt{n}$-consistency and the asymptotic normality in a simple example where a other $\sqrt{n}$-consistent estimator is already known (see [29]).

*Example 3 (Cosines regression function).* Let $f_\theta$ be of the form $f_\theta(x) = \sum_{j=1}^d \theta_j \cos(jx)$ and $p_\varepsilon$ satisfying $(N_2)$ with $\rho \leq 2$. Let $K$ be such $K^*(t) = \mathbb{1}_{|t| \leq 1}$ and let $w(x) = \exp\{-x^2/(4\beta)\}$. Then the conditions $(C_1)$–$(C_3)$ are satisfied and the estimator $\tilde{\theta}_1$ is a $\sqrt{n}$-consistent and asymptotically Gaussian estimator of $\theta^0$.

*Remark 4.2.* This example has already been considered in [29] and [8], but the $\sqrt{n}$-consistency as well as the asymptotic normality is new in this context. More precisely, the estimator constructed in [29] has a rate of convergence of order $\exp(\sqrt{\log n})/n$ for Gaussian errors.

*Example 4 (Cauchy regression function 1).* Consider model (1.1) with $f_\theta(x) = \theta/(1 + x^2)$ satisfying $(R_1)$ with $a = 0, b = 1/2$ and $r = 1$ and $p_\varepsilon$ satisfying $(N_2)$ with $\rho \leq 2$. Let $K$ be such $K^*(t) = \mathbb{1}_{|t| \leq 1}$ and let $w(x) = (1 + x^2)^4 \exp\{-x^2/(4\beta)\}$. With our choice of $w$, the functions $f_\theta w$, $f_\theta^2 w$ and their derivatives in $\theta$ up to order 3 satisfy $(R_1)$ with $\rho < r = 2$ or $\rho = r = 2$ and $b > \beta$. Consequently, the conditions $(C_1)$–$(C_3)$ are satisfied and the estimator $\tilde{\theta}_1$ is a $\sqrt{n}$-consistent and asymptotically Gaussian estimator of $\theta^0$.

This simple example underlies the importance of the smoothing weight function $w$ in the construction of $\tilde{\theta}_1$. Indeed, for a Gaussian noise $\varepsilon$, without a smoothing function $w$ in front of the regression function, Theorem 3.1 predicts (as for the estimator in [29]) a rate of convergence of order $\exp(-\sqrt{\log n})$ instead of the parametric rate of convergence.

**Example 5 (Laplace regression function).** *Consider model (1.1) with $f_\theta(x) = \theta f(x)$ and $f(x) = \exp(-|x|/2)$. The Fourier transform of $f$ and hence of $f_\theta$ is slowly decaying, like $|u|^{-2}$ as $|u| \to \infty$. The estimator $\tilde\theta_1$ with $w \equiv 1$ would not provide $\sqrt n$ consistent estimator for smoother noise densities (as soon as $|p_\varepsilon^*(u)| \le o(|u|^{-2})$ with $|u| \to \infty$). A closer look tells us that $f_\theta$ and its derivative in $\theta$ is $\mathcal{C}^\infty$ except at one point $x = 0$. Therefore, a proper choice of $w$ can smooth out at 0 and make $wf_\theta$, $wf_\theta^2$ and their derivatives in $\theta$ infinitely differentiable functions in $x$. With such a choice of the weight function $w$, the estimator $\tilde\theta_1$ attains the parametric rate of convergence for a large set of noise densities, such that $p_\varepsilon$ satisfies (N$_2$) for some $0 < \rho < 1$. Even if $\rho \ge 1$ and the noise is smoother than that (e.g., Gaussian for $\rho = 2$), the rate of $\tilde\theta_1$ is much faster when using our choice of $w$ then it would be for $w \equiv 1$ or when using the estimator proposed in [29].*

Let us be more precise on the suitable choice of $w$ and define

$$\Psi_{a,b}(x) = \exp\left(-\frac{1}{(x-a)^R(b-x)^R}\right) I_{[a,b]}(x), \tag{4.1}$$

where $-\infty < a < b < \infty$ are fixed and $R > 0$. Following Lepski and Levit [23] and Fedoryuk [11], p. 346, Theorem 7.3, the Fourier transform of this function is such that

$$|\Psi_{a,b}^*(u)| \le c\exp(-C|u|^{R/(R+1)}) \quad \text{as } |u| \to \infty$$

and $c, C > 0$ are constants. Then take $w$ like $\Psi_{0,100}$ or $\Psi_{-100,0}$ or their sum (for a $R > 0$ large enough).

Another way to smooth without restraining to compact support is the following. Let $w(x) = \exp(-1/|x|^{2R})$ a weight function which smoothes at 0 as $R > 0$ is large.

We can vary the coefficient in the exponential in the expression of $w$ and check that $f_\theta w$, $f_\theta^2 w$ and their derivatives up to order 3 satisfy (R$_1$) with the same $r = R/(R+1)$ closer to 1 as $R$ becomes large and $b > 0$.

If the noise satisfies (N$_2$) with $0 \le \rho < 1$, then we find $R$ large enough such that $r = R/(R+1) > \rho$, and thus the conditions (C$_1$)–(C$_3$) are satisfied. Consequently, the estimator $\tilde\theta_1$ is a $\sqrt n$-consistent and asymptotically Gaussian estimator of $\theta^0$.

If $\rho \ge 1$, for this choice of $w$, then

$$\mathbb{E}\|\tilde\theta_1 - \theta^0\|_{\ell^2}^2 = \mathrm{O}(1)(\log n)^{(1-2a-r)/\rho} \exp\{-2b(\log n/(2\beta))^{r/\rho}\}.$$

Note that for the same regression functions with $f(x) = \exp(|x|)$ we can multiply the previous weight function $w$ by $\exp(-4|x|)$ or by $\exp(-x^2)$ in order to solve integrability problems without changing the previous conclusions.

**Example 6 (Irregular regression function).** *Consider model (1.1) with $f_\theta(x) = \theta \mathbb{1}_{[-1,1]}(x)$ and $p_\varepsilon$ satisfying (N$_2$). Let $K$ be such that $K^*(t) = \mathbb{1}_{|t|\le 1}$ and take $w = \Psi_{-1,1}$ for a $R > 0$ defined by (4.1).*

*If $\rho = 0$ in (N$_2$), then the conditions (C$_1$)–(C$_3$) are satisfied. Consequently, the estimator $\tilde\theta_1$ is a $\sqrt n$-consistent and asymptotically Gaussian estimator of $\theta^0$.*

*If $\rho > 0$, then the best rate for estimating $\theta^0$ is obtained by choosing $w = \Psi_{-1,1}$ with $R > 0$ sufficiently large such that $wf_\theta$ and $wf_\theta^2$ satisfy (R$_1$) with $0 < r = R/(R+1) < 1$ as close to 1 as needed.*

*It follows that if $0 < \rho < 1$, then we can find $w = \Psi_{-1,1}$ (with $R$ large enough) satisfying (R$_1$) with $r = R/(R+1) > \rho$, and hence the conditions (C$_4$)–(C$_7$) as well as conditions (C$_1$)–(C$_3$) are satisfied. Thus, $\tilde\theta_1$ is a $\sqrt n$-consistent and asymptotically Gaussian estimator of $\theta^0$. Whereas, if $\rho \ge 1$, for a suitably chosen $w$, then*

$$\mathbb{E}\|\tilde\theta_1 - \theta^0\|_{\ell^2}^2 = \mathrm{O}(1)(\log n)^{(1-2a-r)/\rho} \exp\{-2b(\log n/(2\beta))^{r/\rho}\}.$$

**Example 7 (Polygonal regression function).** *Consider model (1.1) with $f_\theta(x) = \theta_0 + \theta_1 x + \theta_2(x-a)_+ + \theta_3|x-b|^3$ and $p_\varepsilon$ satisfying (N$_2$). Let $K$ be such $K^*(t) = \mathbb{1}_{|t|\le 1}$. This regression function is $\mathcal{C}^\infty$ except at points $a$ and $b$ where it is not differentiable. For $R > 0$, let $w(x) = \Psi_{a-100,a}(x) + \Psi_{a,b}(x) + \Psi_{b,b+100}(x)$, with $\Psi$* *defined in (4.1). The idea is that we can truncate the regression function (say on the interval $[a-100, b+100]$, or larger) in order to smooth out at points $a$, $b$ and the end points of the support of $w$.*

*If the noise satisfies $(N_2)$ with $0 \leq \rho < 1$, then take $R$ large enough such that $r = R/(R+1) > \rho$, and thus the conditions $(C_1)$–$(C_3)$ are satisfied. Consequently, the estimator $\tilde{\theta}_1$ is a $\sqrt{n}$-consistent and asymptotically Gaussian estimator of $\theta^0$.*

*If $\rho \geq 1$ in $(N_2)$, according to Table 1*

$$\mathbb{E}\|\tilde{\theta}_1 - \theta^0\|_{\ell^2}^2 = O(1)(\log n)^{(1-2a-r)/\rho} \exp\{-2b(\log n/(2\beta))^{r/\rho}\}.$$

*Comments on the examples 5, 6 and 7*

In those three examples, $\tilde{\theta}_1$ achieves the $\sqrt{n}$-rate of convergence provided that $p_\varepsilon$ is ordinary smooth or super smooth with an exponent $\rho < 1$. But, $f_\theta w$ will satisfy $(R_1)$ with $r$ at most such that $r < 1$ and it seems, therefore, impossible to have $(w_\theta f_\theta)^*/p_\varepsilon^*$ in $\mathbb{L}_1(\mathbb{R})$ if the $\varepsilon_i$'s are Gaussian. For this regression functions, if the $\varepsilon_i$'s are Gaussian, the least square criterion can not be estimated with the parametric rate of convergence, and hence could probably not provide a $\sqrt{n}$-consistent estimator of $\theta^0$ in this context. Nevertheless, even in cases where the parametric rate of convergence seems not achievable by such an estimator, the resulting rate of the risk of $\tilde{\theta}_1$ is clearly infinitely faster than the logarithmic rate we could have without a proper of $w$ or by using the estimator proposed by Taupin [29].

## 5. Second estimation procedure and comments on the conditions ensuring the $\sqrt{n}$-consistency

### 5.1. Construction and study of the risk of a second estimator

We now propose another estimator of $\theta^0$, based on sufficient conditions allowing to construct directly a $\sqrt{n}$-consistent estimator of $\tilde{S}_{\theta^0,g}$ defined in (2.1), based on $(Y_i, Z_i)$, $i = 1, \ldots, n$.

We say that the conditions $(C_4)$–$(C_7)$ hold if there exists a weight function $w$ and there exist functions $\tilde{\Phi}_{\theta,\varepsilon,j}$, $j = 1, 2, 3$ not depending on $g$, such that for all $\theta \in \Theta$ and for all $g$

$(C_4)$ $\qquad \int y^2 w(x) f_{Y,X}(y,x) \, dy \, dx = \int y^2 \tilde{\Phi}_{\theta,\varepsilon,3}(z) f_{Y,Z}(y,z) \, dy \, dz,$

$\qquad\qquad \int y f_\theta(x) w(x) f_{Y,X}(y,x) \, dy \, dx = \int y \tilde{\Phi}_{\theta,\varepsilon,2}(z) f_{Y,Z} \, dy \, dz$

$\qquad$ and $\quad \int w(x) f_\theta^2(x) g(x) \, dx = \int \tilde{\Phi}_{\theta,\varepsilon,1}(z) h(z) \, dz;$

$(C_5)$ For $j = 1, 2, 3, \mathbb{E}[\sup_{\theta \in \Theta} |\frac{\partial \tilde{\Phi}_{\theta,\varepsilon,j}(Z)}{\partial \theta}|] < \infty$;

$(C_6)$ For $j = 1, 2, 3$ and for all $\theta \in \Theta, \mathbb{E}[|\frac{\partial \tilde{\Phi}_{\theta,\varepsilon,j}}{\partial \theta}(Z)|^2] < \infty$;

$(C_7)$ For $j = 1, 2, 3$ and for all $\theta \in \Theta, \mathbb{E}[|\frac{\partial^2 \tilde{\Phi}_{\theta,\varepsilon,j}}{\partial \theta^2}(Z)|] < \infty$.

Note that $\tilde{\Phi}_{\theta,\varepsilon,3}$ exists as soon as the chosen weight function $w$ is smoother than $p_\varepsilon$ in the way that $w^*/p_\varepsilon^*$ belongs to $\mathbb{L}_1(\mathbb{R})$. Furthermore, $\tilde{\Phi}_{\theta,\varepsilon,3} \equiv \tilde{\Phi}_{\varepsilon,3}$ does not depend on $\theta$. We refer to Section 5.2 for details on how to construct such functions $\tilde{\Phi}_{\theta,\varepsilon,j}$.

Under $(C_4)$–$(C_7)$, we propose to estimate $\tilde{S}_{\theta^0,g}(\theta)$ by

$$\tilde{S}_{n,2}(\theta) = \frac{1}{n} \sum_{i=1}^n [(Y_i^2 \Phi_{\theta,\varepsilon,3}(Z_i) - 2Y_i \Phi_{\theta,\varepsilon,2}(Z_i) + \Phi_{\theta,\varepsilon,1}(Z_i)], \tag{5.2}$$

and hence $\theta^0$ is estimated by

$$\tilde{\theta}_2 = \underset{\theta \in \Theta}{\arg\min}\, \tilde{S}_{n,2}(\theta). \tag{5.3}$$

**Remark 5.1.** *The main difficulty for finding such functions $\tilde{\Phi}_{\theta,\varepsilon,j}$ lies in the constraint that we expect that they do not depend on the unknown density $g$. If we relax this constraint, obviously there are a lot of solutions.*

**Theorem 5.1.** *Consider model* (1.1) *under the assumptions* $(A_1)$–$(A_4)$, $(II_1)$, $(II_2)$, $(N_1)$ *and the conditions* $(C_4)$–$(C_7)$. *Then $\tilde{\theta}_2$, defined by* (5.3) *is a $\sqrt{n}$-consistent estimator of $\theta^0$ which satisfies moreover that $\sqrt{n}(\tilde{\theta}_2 - \theta^0) \xrightarrow[n\to\infty]{\mathcal{L}} \mathcal{N}(0, \tilde{\Sigma}_2)$, with $\tilde{\Sigma}_2$ that equals*

$$\left(\mathbb{E}\left[2w(X)\left(\frac{\partial f_\theta(X)}{\partial \theta}\right)\left(\frac{\partial f_\theta(X)}{\partial \theta}\right)^\top\right]\bigg|_{\theta=\theta^0}\right)^{-1} \tilde{\Sigma}_{0,2} \left(\mathbb{E}\left[2w(X)\left(\frac{\partial f_\theta(X)}{\partial \theta}\right)\left(\frac{\partial f_\theta(X)}{\partial \theta}\right)^\top\right]\bigg|_{\theta=\theta^0}\right)^{-1},$$

*where $\tilde{\Sigma}_{0,2}$ equals*

$$\mathbb{E}\left[\left(\frac{\partial(\tilde{\Phi}_{\theta,\varepsilon,1}(Z) - 2Y\tilde{\Phi}_{\theta,\varepsilon,2}(Z))}{\partial \theta}\bigg|_{\theta=\theta^0}\right)\left(\frac{\partial(\tilde{\Phi}_{\theta,\varepsilon,1}(Z) - 2Y\tilde{\Phi}_{\theta,\varepsilon,2}(Z))}{\partial \theta}\bigg|_{\theta=\theta^0}\right)^\top\right].$$

**Remark 5.2.** *In the Example 2, one can also choose $w \equiv 1$ and use that*

$$\mathbb{E}[\exp(\theta X)] = \frac{\mathbb{E}[\exp(\theta Z)]}{\mathbb{E}[\exp(\theta \varepsilon)]}.$$

*This implies that if we denote by*

$$\tilde{\Phi}_{\theta,\varepsilon,1}(Z) = \frac{\exp(2\theta Z)}{\mathbb{E}[\exp(2\theta\varepsilon)]} \quad \text{and} \quad \tilde{\Phi}_{\theta,\varepsilon,2}(Z) = \frac{\exp(\theta Z)}{\mathbb{E}[\exp(\theta\varepsilon)]}$$

*then $\mathbb{E}_h[\tilde{\Phi}_{\theta,\varepsilon,1}(Z)] = \mathbb{E}_g[f_\theta^2(X)]$ and $\mathbb{E}_{\theta^0,h}[Y\tilde{\Phi}_{\theta,\varepsilon,2}(Z)] = \mathbb{E}_{\theta^0,g}[Y f_\theta(X)]$. Consequently, $\tilde{S}_{n,2}$ satisfies*

$$\tilde{S}_{n,2}(\theta) = \frac{1}{n}\sum_{i=1}^n \left[Y_i - \frac{\exp(\theta Z_i)}{\mathbb{E}[\exp(\theta\varepsilon)]}\right]^2 = \frac{1}{n}\sum_{i=1}^n [Y_i^2 - 2Y_i\tilde{\Phi}_{\theta,\varepsilon,2}(Z_i) + \tilde{\Phi}_{\theta,\varepsilon,1}(Z_i)]. \tag{5.4}$$

*In this case $\tilde{\theta}_2$ is also a $\sqrt{n}$-consistent and asymptotically Gaussian estimator of $\theta^0$.*

**Remark 5.3.** *In the Example 3, one can also choose $w \equiv 1$ and use that $\mathbb{E}[\exp(ijX)] = \mathbb{E}[\exp(ijZ)]/\mathbb{E}[\exp(ij\varepsilon)]$. This implies that if we denote by*

$$\tilde{\Phi}_{\theta,\varepsilon,1}(Z) = \frac{1}{4}\Bigg\{1 + \sum_{j=1}^d \theta_j^2\left[\frac{\exp(2ijZ)}{p_\varepsilon^*(2j)} + \frac{\exp(-2ijZ)}{p_\varepsilon^*(-2j)}\right]$$

$$+ \sum_{j=1}^d \sum_{k\neq j} \theta_j\theta_k\left[\frac{\exp(i(j+k)Z)}{p_\varepsilon^*(j+k)} + \frac{\exp(-i(j+k)Z)}{p_\varepsilon^*(-(j+k))} + \frac{\exp(i(j-k)Z)}{p_\varepsilon^*(j-k)} + \frac{\exp(i(-j+k)Z)}{p_\varepsilon^*(-j+k)}\right]\Bigg\}$$

*and*

$$\tilde{\Phi}_{\theta,\varepsilon,2}(Z) = \frac{1}{2}\left[\frac{\exp(ijZ)}{p_\varepsilon^*(j)} + \frac{\exp(-ijZ)}{p_\varepsilon^*(-j)}\right]$$

then $\mathbb{E}[\tilde{\Phi}_{\theta,\varepsilon,1}(Z)] = \mathbb{E}[f_\theta^2(X)]$ and $\mathbb{E}[Y\tilde{\Phi}_{\theta,\varepsilon,2}(Z)] = \mathbb{E}[Yf_\theta(X)]$. Consequently, $\tilde{S}_{n,2}$ satisfies

$$\tilde{S}_{n,2}(\theta) = \frac{1}{n}\sum_{i=1}^{n}[Y_i^2 - 2Y_i\tilde{\Phi}_{\theta,\varepsilon,2}(Z_i) + \tilde{\Phi}_{\theta,\varepsilon,1}(Z_i)]. \tag{5.5}$$

In this case $\tilde{\theta}_2$ with $w \equiv 1$ is again a $\sqrt{n}$-consistent and asymptotically Gaussian estimator of $\theta^0$. In the same way, $\tilde{\theta}_1$ with $w \equiv 1$ is again a $\sqrt{n}$-consistent and asymptotically Gaussian estimator of $\theta^0$.

**Remark 5.4.** Consider the Example 4. If we choose $w$ as for $\tilde{\theta}_1$, then the conditions (C$_4$)–(C$_7$) are fulfilled and $\tilde{\theta}_2$ is a $\sqrt{n}$-consistent and asymptotically Gaussian estimator of $\theta^0$. Consequently, in this example, the conditions (5.1) given in [9] and [10] are not satisfied when our conditions (C$_4$)–(C$_7$) hold. Hence our estimation procedure provides a $\sqrt{n}$-consistent estimator when the estimation procedure in [10] and [9] fails.

### 5.2. Comments on the conditions ensuring the $\sqrt{n}$-consistency

**Comment 1.** Let us briefly compare the conditions (C$_1$)–(C$_3$) to the conditions (C$_4$)–(C$_7$). It is noteworthy that the conditions (C$_4$)–(C$_7$) are more general. For instance condition (C$_1$) implies (C$_4$), with $\tilde{\Phi}_{\theta,\varepsilon,j}$ defined by $\tilde{\Phi}^*_{\theta,\varepsilon,1} = (wf_\theta^2)^*/p_\varepsilon^*$, $\tilde{\Phi}^*_{\theta,\varepsilon,2} = (wf_\theta)^*/p_\varepsilon^*$ and $\tilde{\Phi}^*_{\theta,\varepsilon,3} = w^*/p_\varepsilon^*$. This comes from the following equalities $\mathbb{E}[\tilde{\Phi}_{\theta,\varepsilon,1}(Z)] = \mathbb{E}[(wf_\theta^2)(X)]$,

$$\begin{aligned}\mathbb{E}[Y^2\tilde{\Phi}_{\theta,\varepsilon,3}(Z)] &= \mathbb{E}[f_{\theta^0}^2(X)\tilde{\Phi}_{\theta,\varepsilon,3}(Z)] + \sigma_{\xi,2}^2\mathbb{E}[\tilde{\Phi}_{\theta,\varepsilon,3}(Z)]\\ &= \langle f_{\theta_0}^2 g, \tilde{\Phi}_{\theta,\varepsilon,3}\star p_\varepsilon\rangle + \sigma_{\xi,2}^2\langle g, \tilde{\Phi}_{\theta,\varepsilon,3}\star p_\varepsilon\rangle\\ &= (2\pi)^{-1}\langle (f_{\theta_0}^2 g)^*, (\tilde{\Phi}_{\theta,\varepsilon,3})^*p_\varepsilon^*\rangle + \sigma_{\xi,2}^2(2\pi)^{-1}\langle g^*, (\tilde{\Phi}_{\theta,\varepsilon,3})^*p_\varepsilon^*\rangle\\ &= \langle f_{\theta^0}^2 g, w\rangle + \sigma_{\xi,2}^2\langle g, w\rangle = \mathbb{E}[Y^2 w(X)]\end{aligned}$$

and

$$\begin{aligned}\mathbb{E}[Y\tilde{\Phi}_{\theta,\varepsilon,2}(Z)] &= \langle f_{\theta_0}g, \tilde{\Phi}_{\theta,\varepsilon,2}\star p_\varepsilon\rangle = (2\pi)^{-1}\langle (f_{\theta_0}g)^*, (\tilde{\Phi}_{\theta,\varepsilon,2})^*p_\varepsilon^*\rangle\\ &= \langle f_{\theta^0}g, f_\theta w\rangle = \mathbb{E}[Y(wf_\theta)(X)].\end{aligned}$$

But in the condition (C$_4$), $f_\theta w$, $f_\theta^2 w$ and $w$ are not necessarily in $\mathbb{L}_1(\mathbb{R})$. Consequently, under (C$_1$)–(C$_3$), the conditions (C$_1$)–(C$_4$) hold and, by denoting $\tilde{\Phi}_{\theta,\varepsilon,1} = (wf_\theta^2)^*/p_\varepsilon^*$, $\tilde{\Phi}_{\theta,\varepsilon,2} = (wf_\theta)^*/p_\varepsilon^*$ and $\tilde{\Phi}_{\theta,\varepsilon,3} = w^*/p_\varepsilon^*$, we get that $\tilde{\Sigma}_{0,1} = \tilde{\Sigma}_{0,2}$ with $\tilde{\Sigma}_{0,1}$ and $\tilde{\Sigma}_{0,2}$ defined in Theorems 3.2 and 5.1. Nevertheless, the conditions (C$_1$)–(C$_3$) are more tractable and constructive conditions.

**Comment 2.** The conditions (C$_4$)–(C$_7$) have to be related to the conditions given in [9, 10] and [1] in the context of the functional errors in variables model with $X_1,\ldots,X_n$ not random. In both mentionned papers, in order to construct $\sqrt{n}$-consistent estimator, they assume that there exist two functions $\phi_1$ and $\phi_2$ such that

$$\mathbb{E}(\phi_1(x+\varepsilon,\theta)) = f_\theta(x) \quad \text{and} \quad \mathbb{E}(\phi_2(x+\varepsilon,\theta)) = f_\theta^2(x). \tag{5.1}$$

Clearly our conditions (C$_4$)–(C$_7$) are less restrictive than the condition (5.1), and hence they allow to achieve the parametric rate of convergence for various type of regression functions, by using the possibility of the choice of the weight function $w_\theta$.

For instance, let us reconsider Example 4 where $f_\theta(x) = \theta/(1+x^2)$ and $p_\varepsilon$ is the Gaussian density. Then the conditions given in [9] and [10] are not fulfilled. Whereas $\mathbb{E}[f_\theta(X)w(X)]$ and $\mathbb{E}[f_\theta^2(X)w(X)]$ can be estimated with the parametric rate of convergence, by taking $w(x) = (1+x^2)^4\exp(-x^2/(4\beta))$. It follows from the fact that condition (C$_4$)–(C$_7$) are fulfilled in this special example (see previous comment for the construction of auxiliary functions $\tilde{\Phi}_{\theta,\varepsilon,j}$, $j=1$ to 3 appearing in the condition). Nevertheless, such weight function are not always available, and, therefore, those conditions (C$_4$)–(C$_7$) are not always fulfilled.

## 6. Extensions of previous estimation procedures

In the estimation procedure previously presented, the weight function is used in order to make $wf_\theta$, $wf_\theta^2$ integrable and also in order to get smooth $wf_\theta$, $wf_\theta^2$ and derivatives in $\theta$. In a large class of regression functions, the weight function can smooth these functions without depending on $\theta$. Sometimes, the smoothing properties and hence the rate of convergence are improved by making $w$ to depend on $\theta$. It appears, in particular, when the points where we need to smooth are related to $\theta$.

The second estimation procedure uses an estimator based on the observations $(Y_i, Z_i)$, $i = 1, \ldots, n$, of the least square contrast

$$S_{\theta^0,g}(\theta) = \mathbb{E}[((Y - f_\theta(X))^2 - \sigma_{\xi,2}^2) w_\theta(X)], \tag{6.1}$$

where $w_\theta$ is a positive weight function to be suitably chosen. This criterion, which requires the knowledge of $\mathrm{Var}(\xi)$, actually writes $\mathbb{E}[((f_{\theta^0}(X) - f_\theta(X))^2 w_\theta(X)]$ and it is minimized in $\theta = \theta^0$ for any positive weight function $w_\theta$. Subsequently, we consider the following assumptions.

*Identifiability and moment assumptions*

(I2$_1$) The variance $\sigma_{\xi,2}^2 = \mathrm{Var}(\xi_i)$ is known.
(I2$_2$) The quantity $S_{\theta^0,g}(\theta) = \mathbb{E}[w_\theta(X)(Y - f_\theta(X))^2] - \sigma_{\xi,2}^2 \mathbb{E}(w_\theta(X))$ admits one unique minimum at $\theta = \theta^0$.
(I2$_3$) For all $\theta \in \Theta$ the matrix $S_{\theta^0,g}^{(2)}(\theta) = \partial^2 S_{\theta^0,g}(\theta)/\partial \theta^2$ exists and

$$S_{\theta^0,g}^{(2)}(\theta^0) = \mathbb{E}\left[2 w_{\theta^0}(X) \left(\frac{\partial f_\theta(X)}{\partial \theta}\bigg|_{\theta=\theta^0}\right) \left(\frac{\partial f_\theta(X)}{\partial \theta}\bigg|_{\theta=\theta^0}\right)^\top \right] \quad \text{is positive definite.}$$

(A$_5$) The quantities $\mathbb{E}[w_\theta^2(X)(Y - f_\theta(X))^4]$ and their derivatives up to order 2 with respect to $\theta$ are finite.

According to (2.2), $S_{\theta^0,g}(\theta)$ is naturally estimated by the empirical criterion

$$S_{n,1}(\theta) = \frac{1}{n} \sum_{i=1}^{n} \int [(Y_i - f_\theta(x))^2 - \sigma_{\xi,2}^2] w_\theta(x) K_{n,C_n}(C_n(x - Z_i)) \, dx, \tag{6.2}$$

where $K_{n,C_n}(\cdot) = C_n K_n(C_n \cdot)$ is a deconvolution kernel satisfying (2.4). Using this empirical criterion, under (N$_1$), we propose to estimate $\theta^0$ by

$$\widehat{\theta}_1 = \arg\min_{\theta \in \Theta} S_{n,1}(\theta). \tag{6.3}$$

**Theorem 6.1.** *Let $\widehat{\theta}_1 = \widehat{\theta}_1(C_n)$ be defined by (6.3) under the assumptions (I2$_1$)–(I2$_3$), (N$_1$), (A$_1$)–(A$_4$) and (A$_5$) with $w$ replaced by $w_\theta$.*

(1) *Then for all of the sequences $C_n$ such that for all $\theta \in \Theta$, $f_\theta w_\theta$ and $f_\theta^2 w_\theta$ satisfy (3.6), $\mathbb{E}(\|\widehat{\theta}_1(C_n) - \theta^0\|_{\ell^2}^2) = o(1)$, as $n \to \infty$ and $\widehat{\theta}_1(C_n)$ is a consistent estimator of $\theta^0$.*
(2) *Assume moreover that for all $\theta \in \Theta$, $f_\theta w_\theta$ and $f_\theta^2 w_\theta$ and their derivatives up to order 3 with respect to $\theta$ satisfy (3.6). Then $\mathbb{E}(\|\widehat{\theta}_1 - \theta^0\|_{\ell^2}^2) = O(\varphi_n^2)$ with $\varphi_n$ given $\varphi_n = \|(\varphi_{n,j})\|_{\ell^2}$ with $\varphi_{n,j}^2 = B_{n,j}^2(\theta^0) + V_{n,j}(\theta^0)/n$, $j = 1, \ldots, d$, where*

$$B_{n,j}(\theta) = \min\{B_{n,j}^{[1]}(\theta), B_{n,j}^{[2]}(\theta)\} \quad \text{and} \quad V_{n,j}(\theta) = \min\{V_{n,j}^{[1]}(\theta), V_{n,j}^{[2]}(\theta)\},$$

$$B_{n,j}^{[q]} = \left\|\left(\frac{\partial(w_\theta)}{\partial \theta_j}\right)^*(K_{C_n}^* - 1)\right\|_q^2 + \left\|\left(\frac{\partial(f_\theta w_\theta)}{\partial \theta_j}\right)^*(K_{C_n}^* - 1)\right\|_q^2 + \left\|\left(\frac{\partial(f_\theta^2 w_\theta)}{\partial \theta_j}\right)^*(K_{C_n}^* - 1)\right\|_q^2$$

*and*

$$V_{n,j}^{[q]}(\theta) = \left\|\left(\frac{\partial(w_\theta)}{\partial\theta_j}\right)^* \frac{K_{C_n}^*}{p_\varepsilon^*}\right\|_q^2 + \left\|\left(\frac{\partial(f_\theta w_\theta)}{\partial\theta_j}\right)^* \frac{K_{C_n}^*}{p_\varepsilon^*}\right\|_q^2 + \left\|\left(\frac{\partial(f_\theta^2 w_\theta)}{\partial\theta_j}\right)^* \frac{K_{C_n}^*}{p_\varepsilon^*}\right\|_q^2.$$

**Remark 6.1.** *The Remark 3.2 is still valid.*

*6.1. Consequence: a sufficient condition to obtain the parametric rate of convergence with $\widehat{\theta}_1$*

We say that the conditions $(C_8)$–$(C_{10})$ hold if there exists a weight function $w_\theta$ such that for all $\theta \in \Theta$,

($C_8$) the functions $(w_\theta f_\theta), (w_\theta)$ and $(w_\theta f_\theta^2)$ belong to $\mathbb{L}_1(\mathbb{R})$ and the functions $(w_\theta)^*/p_\varepsilon^*$,
$\sup_\theta (f_\theta w_\theta)^*/p_\varepsilon^*, \sup_\theta (f_\theta^2 w_\theta)^*/p_\varepsilon^*$ belong to $\mathbb{L}_1(\mathbb{R}) \cap \mathbb{L}_2(\mathbb{R})$;

($C_9$) the functions $\sup_{\theta\in\Theta}(\frac{\partial w_\theta}{\partial\theta})^*/p_\varepsilon^*, \sup_{\theta\in\Theta}(\frac{\partial(f_\theta w_\theta)}{\partial\theta})^*/p_\varepsilon^*$ and $\sup_{\theta\in\Theta}(\frac{\partial(f_\theta^2 w_\theta)}{\partial\theta})^*/p_\varepsilon^*$ belong to $\mathbb{L}_1(\mathbb{R}) \cap \mathbb{L}_2(\mathbb{R})$;

($C_{10}$) the functions $(\frac{\partial^2 w_\theta}{\partial\theta^2})^*/p_\varepsilon^*, (\frac{\partial^2(f_\theta w_\theta)}{\partial\theta^2})^*/p_\varepsilon^*$ and $(\frac{\partial^2(f_\theta^2 w_\theta)}{\partial\theta^2})^*/p_\varepsilon^*$ belong to $\mathbb{L}_1(\mathbb{R}) \cap \mathbb{L}_2(\mathbb{R})$.

**Theorem 6.2.** *Consider model* (1.1) *under the assumptions* $(I2_1)$–$(I2_3)$, $(N_1)$, $(A_1)$–$(A_4)$ *for $w$ replaced by $w_\theta$ and under* $(A_5)$, $(C_8)$–$(C_{10})$. *Then $\widehat{\theta}_1$ defined by* (6.3) *is a $\sqrt{n}$-consistent estimator of $\theta^0$ which satisfies moreover that $\sqrt{n}(\widehat{\theta}_1 - \theta^0) \xrightarrow[n\to\infty]{\mathcal{L}} \mathcal{N}(0, \Sigma_1)$, with $\Sigma_1$ that equals*

$$\left(\mathbb{E}\left[2w_\theta(X)\left(\frac{\partial f_\theta(X)}{\partial\theta}\right)\left(\frac{\partial f_\theta(X)}{\partial\theta}\right)^\top\right]\bigg|_{\theta=\theta^0}\right)^{-1} \Sigma_{0,1} \left(\mathbb{E}\left[2w_\theta(X)\left(\frac{\partial f_\theta(X)}{\partial\theta}\right)\left(\frac{\partial f_\theta(X)}{\partial\theta}\right)^\top\right]\bigg|_{\theta=\theta^0}\right)^{-1},$$

*where*

$$\Sigma_{0,1} = (2\pi)^{-2} \mathbb{E}\Bigg\{\left[\int \left(\frac{\partial[f_\theta^2 w_\theta - 2Y f_\theta w_\theta + (Y^2 - \sigma_{\xi,2}^2)w_\theta]}{\partial\theta}\bigg|_{\theta=\theta^0}\right)^*(u) \frac{e^{-iuZ}}{p_\varepsilon^*(u)}\,du\right]$$
$$\times \left[\int \left(\frac{\partial[f_\theta^2 w_\theta - 2Y f_\theta w_\theta + (Y^2 - \sigma_{\xi,2}^2)w_\theta]}{\partial\theta}\bigg|_{\theta=\theta^0}\right)^*(u) \frac{e^{-iuZ}}{p_\varepsilon^*(u)}\,du\right]^\top\Bigg\}.$$

The resulting rate when $(f_\theta w_\theta)$ and $(f_\theta^2 w_\theta)$, as well as their derivatives with respect to $\theta$ up to order 3 satisfy $(R_1)$, are given in the following corollary.

**Corollary 6.1.** *Under the assumptions of Theorem* 6.1, *assume that $p_\varepsilon$ satisfies* $(N_2)$ *and that for all $\theta \in \Theta$, $(f_\theta w_\theta)$ and $(f_\theta^2 w_\theta)$ and their derivatives with respect to $\theta$ up to order 3, satisfy* $(R_1)$. *Then Corollary* 3.1 *still holds with $wf_\theta$ replaced by $w_\theta f_\theta$.*

*6.2. Construction and study of the risk of the estimator $\widehat{\theta}_2$*

In the same spirit as $\widehat{\theta}_1$, we propose another estimator whose construction is based on sufficient conditions allowing to construct a direct and $\sqrt{n}$-consistent estimator of $S_{\theta^0,g}(\theta)$ defined in (6.1), using the observations $(Y_i, Z_i)$, $i = 1, \ldots, n$.

We say that the conditions $(C_{11})$–$(C_{14})$ hold if there exists a weight function $w_\theta$ and there exist three functions $\Phi_{\theta,\varepsilon,j}$, $j = 1, \ldots, 3$, not depending on $g$ such that for all $\theta \in \Theta$ and for all $g$

($C_{11}$) $\quad \int w_\theta(x)g(x)\,dx = \int \Phi_{\theta,\varepsilon,3}(z)h(z)\,dz,$

$$\int y f_\theta(x) w_\theta(x) f_{Y,X}(y,x)\,dy\,dx = \int y \Phi_{\theta,\varepsilon,2}(z) f_{Y,Z}(y,z)\,dy\,dz$$

and $\quad \int w_\theta(x) f_\theta^2(x) g(x)\,dx = \int \Phi_{\theta,\varepsilon,1}(z)h(z)\,dz;$

(C$_{12}$)  For $j = 1, 2, 3$, $\mathbb{E}[\sup_{\theta \in \Theta} |\frac{\partial \Phi_{\theta,\varepsilon,j}}{\partial \theta}(Z)|] < \infty$;

(C$_{13}$)  For $j = 1, 2, 3$, and for all $\theta \in \Theta$, $\mathbb{E}[|\frac{\partial \Phi_{\theta,\varepsilon,j}}{\partial \theta}(Z)|^2] < \infty$;

(C$_{14}$)  For $j = 1, 2, 3$, and for all $\theta \in \Theta$, $\mathbb{E}[|\frac{\partial^2 \Phi_{\theta,\varepsilon,j}}{\partial \theta^2}(Z)|] < \infty$.

Under (C$_{11}$)–(C$_{14}$), we propose to estimate $S_{\theta^0,g}(\theta)$ by

$$S_{n,2}(\theta) = \frac{1}{n} \sum_{i=1}^n [(Y_i^2 - \sigma_{\xi,2}^2)\Phi_{\theta,\varepsilon,3}(Z_i) - 2Y_i \Phi_{\theta,\varepsilon,2}(Z_i) + \Phi_{\theta,\varepsilon,1}(Z_i)]. \tag{6.4}$$

Using this empirical criterion we propose to estimate $\theta^0$ by

$$\widehat{\theta}_2 = \arg\min_{\theta \in \Theta} S_{n,2}(\theta). \tag{6.5}$$

We refer to Section 5.2 for details on how to construct such functions $\Phi_{\theta,\varepsilon,j}$, $j = 1, \ldots, 3$.

***Remark 6.2.*** *As in the first estimation procedure (see Remark 5.1), the main difficulty for finding such function $\Phi_{\theta,\varepsilon,j}$ lies in the constraint that we expect that they do not depend on the unknown density $g$. If we relax this constraint, obviously there are a lot of solutions.*

### 6.3. Asymptotic properties of the estimator $\widehat{\theta}_2$

***Theorem 6.3.*** *Consider the model (1.1) under the assumptions (A$_1$)–(A$_5$), (I2$_1$)–(I2$_3$), (N$_1$) and the conditions (C$_{11}$)–(C$_{14}$). Then $\widehat{\theta}_2$, defined by (6.5) is a $\sqrt{n}$-consistent estimator of $\theta^0$ which satisfies moreover that $\sqrt{n}(\widehat{\theta}_2 - \theta^0) \xrightarrow[n\to\infty]{\mathcal{L}} \mathcal{N}(0, \Sigma_2)$, with $\Sigma_2$ that equals*

$$\left(\mathbb{E}\left[2w_\theta(X)\left(\frac{\partial f_\theta(X)}{\partial \theta}\right)\left(\frac{\partial f_\theta(X)}{\partial \theta}\right)^\top\right]\bigg|_{\theta=\theta^0}\right)^{-1} \Sigma_{0,2} \left(\mathbb{E}\left[2w_\theta(X)\left(\frac{\partial f_\theta(X)}{\partial \theta}\right)\left(\frac{\partial f_\theta(X)}{\partial \theta}\right)^\top\right]\bigg|_{\theta=\theta^0}\right)^{-1},$$

*where $\Sigma_{0,2}$ equals*

$$\mathbb{E}\left[\left(\frac{\partial(\Phi_{\theta,\varepsilon,1}(Z) - 2Y\Phi_{\theta,\varepsilon,2}(Z) + (Y^2 - \sigma_{\xi,2}^2)\Phi_{\theta^0,\varepsilon,3}(Z))}{\partial \theta}\bigg|_{\theta=\theta^0}\right) \times \left(\frac{\partial(\Phi_{\theta,\varepsilon,1}(Z) - 2Y\Phi_{\theta,\varepsilon,2}(Z) + (Y^2 - \sigma_{\xi,2}^2)\Phi_{\theta^0,\varepsilon,3}(Z))}{\partial \theta}\bigg|_{\theta=\theta^0}\right)^\top\right].$$

***Remark 6.3.*** *Also note that by denoting $\Phi_{\theta,\varepsilon,1} = (w_\theta f_\theta^2)^*/p_\varepsilon^*$, $\Phi_{\theta,\varepsilon,2} = (w_\theta f_\theta)^*/p_\varepsilon^*$ and $\Phi_{\theta,\varepsilon,3} = w_\theta^*/p_\varepsilon^*$, we get that $\Sigma_{0,1} = \Sigma_{0,2}$ with $\Sigma_{0,1}$ defined in Theorem 6.2.*

***Remark 6.4.*** *The same comparison between (C$_{11}$)–(C$_{14}$) and (C$_8$)–(C$_{10}$) holds with $w$ replaced by $w_\theta$ and with the $\Phi_{\theta,\varepsilon,j}$ defined by $\Phi_{\theta,\varepsilon,1}^* = (w_\theta f_\theta^2)^*/p_\varepsilon^*$, $\Phi_{\theta,\varepsilon,2}^* = (w_\theta f_\theta)^*/p_\varepsilon^*$ and $\Phi_{\theta,\varepsilon,3}^* = (f_\theta w_\theta)^*/p_\varepsilon^*$ (See Section 5.2).*

## 7. Examples and methodological advice (2)

***Example 8 (Growth curves 1).*** *Consider the Model (1.1) with known $\sigma_{\xi,2}^2$, $f_\theta(x) = \theta_1/(1 + \theta_2 \exp(\theta_3 x))$ and $p_\varepsilon$ satisfying (N$_2$) with $\rho \leq 2$. Let $K^*(t) = \mathbb{1}_{|t| \leq 1}$ and let $w_\theta(x) = (1 + \theta_2 \exp(\theta_3 x))^4 \exp\{-x^2/(4\beta)\}$. Then the conditions (C$_8$)–(C$_{10}$) as well as conditions (C$_{11}$)–(C$_{14}$) are satisfied. Consequently, the estimators $\widehat{\theta}_1$ and $\widehat{\theta}_2$ are $\sqrt{n}$-consistent and asymptotically Gaussian estimators of $\theta^0$, with the same asymptotic variance.*

**Example 9 (Growth curve 2).** *Consider the Model (1.1) with known $\sigma^2_{\xi,2}$, $f_\theta(x) = \theta_2 + (\theta_1 - \theta_2)/(1 + \exp(\theta_3 + \theta_4 x))$ and $p_\varepsilon$ satisfying (N$_2$) with $\rho \leq 2$. Let $K$ be such $K^*(t) = \mathbb{1}_{|t|\leq 1}$ and let $w_\theta(x) = (1 + \exp(\theta_3 + \theta_4 x))^4 \exp\{-x^2/(4\beta)\}$. Then the conditions (C$_8$)–(C$_{10}$) as well as conditions (C$_{11}$)–(C$_{14}$) are satisfied. Consequently, the estimators $\widehat{\theta}_1$ and $\widehat{\theta}_2$ are $\sqrt{n}$-consistent and asymptotically Gaussian estimators of $\theta^0$, with the same asymptotic variance.*

**Example 10 (Cauchy regression function 2).** *Consider the model (1.1) with known $\sigma^2_{\xi,2}$, $f_\theta(x) = 1/(1 + \theta x^2)$ and $p_\varepsilon$ satisfying (N$_2$) with $\rho \leq 2$. Let $K$ be such $K^*(t) = \mathbb{1}_{|t|\leq 1}$ and let $w_\theta(x) = (1 + \theta x^2)^2 \exp\{-x^2/(4\beta)\}$. Then the conditions (C$_8$)–(C$_{10}$) as well as conditions (C$_{11}$)–(C$_{14}$) are satisfied. Consequently, the estimators $\widehat{\theta}_1$ and $\widehat{\theta}_2$ are $\sqrt{n}$-consistent and asymptotically Gaussian estimators of $\theta^0$, with the same asymptotic variance.*

In these three last examples, we see the importance of the weight function $w_\theta$, with the improvement of the rate of convergence by taking $w_\theta$ depending on $\theta$. Clearly, for such regression function, the estimator in [29] does not achieve the parametric rate of convergence whereas, $\widehat{\theta}_1$ or $\widehat{\theta}_2$ do.

## 8. Proofs of theorems

We only detail the proof of Theorems 6.1, 6.2 and 6.3 as well as Corollary 6.1. The proofs of Theorems 3.1, 3.2, and 5.1 and Corollary 3.1 follow the same lines with $w_\theta$ replaced by $w$ and by putting $\sigma_{\xi,2} = 0$ in the proofs (and only in the proofs). In the sequel, we denote by $C$ an absolute constant whose value may vary from one line to the other, and we always mention the dependency of a constant $C$ with respect to parameters. For instance $C(a,b)$ stands for a constant depending on $a$ and $b$.

### 8.1. Proof of (1) of Theorem 6.1

The main point of the proof consists in showing that for any $\theta$ in $\Theta$, $\mathbb{E}[(S_{n,1}(\theta) - S_{\theta^0,g}(\theta))^2] = \mathrm{o}(1)$, with $S_{\theta^0,g}(\theta)$ admitting a unique minimum in $\theta = \theta^0$. The second part of the proof consists in studying $\omega_2(n,\rho)$ defined as

$$\omega_2(n,\rho) = \sup\{|S_{n,1}(\theta) - S_{n,1}(\theta')| \colon \|\theta - \theta'\|_{\ell^2} \leq \rho\}.$$

By using the regularity assumptions on the regression function $f$, we state that there exist two sequences $\rho_k$ and $\epsilon_k$ tending to 0, such that for all $k \in \mathbb{N}$

$$\lim_{n\to\infty} \mathbb{P}[\omega_2(n,\rho_k) > \epsilon_k] = 0 \quad \text{and that} \quad \mathbb{E}[(\omega_2(n,\rho_k))^2] = \mathrm{O}(\rho_k^2). \tag{8.1}$$

Let us start with the proof of

$$\mathbb{E}[(S_{n,1}(\theta) - S_{\theta^0,g}(\theta))^2] = \mathrm{o}(1), \quad \text{for any } \theta \in \Theta \text{ and as } n \to \infty. \tag{8.2}$$

We have to check successively that for any $\theta \in \Theta$,

$$\mathbb{E}[S_{n,1}(\theta)] - S_{\theta^0,g}(\theta) = \mathrm{o}(1) \quad \text{and} \quad \mathrm{Var}(S_{n,1}(\theta)) = \mathrm{o}(1), \quad \text{as } n \to \infty. \tag{8.3}$$

For both the bias and the variance, we give two upper bounds, based on the two following applications of the Hölder's inequality

$$|\langle \varphi_1, \varphi_2 \rangle| \leq \|\varphi_1\|_2 \|\varphi_2\|_2 \quad \text{and} \quad |\langle \varphi_1, \varphi_2 \rangle| \leq \|\varphi_1\|_\infty \|\varphi_2\|_1. \tag{8.4}$$

**Proof of the first part of (8.3).** According to Lemma A.2, we have

$$\mathbb{E}[S_{n,1}(\theta)] = \mathbb{E}[((Y - f_\theta)^2 w_\theta - \sigma^2_{\xi,2} w_\theta) \star K_{n,C_n}(Z)]$$
$$= \mathbb{E}[(Y^2 - \sigma^2_{\xi,2}) w_\theta \star K_{C_n}(X)] - 2\mathbb{E}[Y(f_\theta w_\theta) \star K_{C_n}(X)] + \mathbb{E}[(f_\theta^2 w_\theta) \star K_{C_n}(X)],$$

and hence

$$\mathbb{E}[S_{n,1}(\theta)] - S_{\theta^0,g}(\theta) = (2\pi)^{-1}\langle (f_{\theta^0}^2 g)^*, w_\theta^*(K_{C_n}^*(\cdot)-1)\rangle - \pi^{-1}\langle (f_{\theta^0}g)^*, (f_\theta w_\theta)^*(K_{C_n}^*-1)\rangle$$
$$+ (2\pi)^{-1}\langle g^*, (f_\theta^2 w_\theta)^*(K_{C_n}^*-1)\rangle.$$

Consequently, under the assumption (A$_3$), $|\mathbb{E}[S_{n,1}(\theta)] - S_{\theta^0,g}(\theta)|^2 = o(1)$. Indeed, according to (8.4), a first upper bound of $|\mathbb{E}[S_{n,1}(\theta)] - S_{\theta^0,g}(\theta)|$ is given by

$$(2\pi)^{-1}\|(f_{\theta^0}^2 g)^*\|_2 \|w_\theta^*(K_{C_n}^*-1)\|_2 + \pi^{-1}\|(f_{\theta^0}g)^*\|_2\|(f_\theta w_\theta)^*(K_{C_n}^*-1)\|_2$$
$$+ (2\pi)^{-1}\|g^*\|_2\|(f_\theta^2 w_\theta)^*(K_{C_n}^*-1)\|_2$$

and hence

$$|\mathbb{E}[S_{n,1}(\theta)] - S_{\theta^0,g}(\theta)|^2$$
$$\leq C(f_{\theta^0}, w_\theta)[\|w_\theta^*(K_{C_n}^*-1)\|_2^2 + \|(w_\theta f_\theta)^*(K_{C_n}^*-1)\|_2^2 + \|(w_\theta f_\theta^2)^*(K_{C_n}^*-1)\|_2^2]. \tag{8.5}$$

And according to the second part of (8.4) we also have that $|\mathbb{E}[S_{n,1}(\theta)] - S_{\theta^0,g}(\theta)|$ is bounded by

$$(2\pi)^{-1}\|(f_{\theta^0}^2 g)^*\|_\infty \|w_\theta^*(K_{C_n}^*-1)\|_1 + \pi^{-1}\|(f_{\theta^0}g)^*\|_\infty\|(f_\theta w_\theta)^*(K_{C_n}^*-1)\|_1$$
$$+ (2\pi)^{-1}\|g^*\|_\infty\|(f_\theta^2 w_\theta)^*(K_{C_n}^*-1)\|_1$$

and hence

$$|\mathbb{E}[S_{n,1}(\theta)] - S_{\theta^0,g}(\theta)|^2$$
$$\leq C(f_{\theta^0}, w_\theta)[\|w_\theta^*(K_{C_n}^*-1)\|_1^2 + \|(w_\theta f_\theta)^*(K_{C_n}^*-1)\|_1^2 + \|(w_\theta f_\theta^2)^*(K_{C_n}^*-1)\|_1^2]. \tag{8.6}$$

Consequently, by combining the two bounds (8.5) and (8.6), we get that $|\mathbb{E}[S_{n,1}(\theta)] - S_{\theta^0,g}(\theta)|^2$ is bounded by

$$C(f_{\theta^0}, w_\theta) \min\{\|w_\theta^*(K_{C_n}^*-1)\|_2^2 + \|(w_\theta f_\theta)^*(K_{C_n}^*-1)\|_2^2 + \|(w_\theta f_\theta^2)^*(K_{C_n}^*-1)\|_2^2,$$
$$\|w_\theta^*(K_{C_n}^*-1)\|_1^2 + \|(w_\theta f_\theta)^*(K_{C_n}^*-1)\|_1^2 + \|(w_\theta f_\theta^2)^*(K_{C_n}^*-1)\|_1^2\}$$

that is by applying Lemma A.1,

$$|\mathbb{E}[S_{n,1}(\theta)] - S_{\theta^0,g}(\theta)|^2 \leq C(f_{\theta^0}, w_\theta, b, r) C_n^{-2a+(1-r)_+(1-r)_-} e^{-2bC_n^r},$$

and the first part of (8.3) follows. $\square$

**Proof of the second part of (8.3).** Since the variables are i.i.d. random variables, the stochastic term on the left-hand side in (8.3) has variance

$$\mathrm{Var}[S_{n,1}(\theta)] \leq n^{-1}\mathbb{E}[|((Y-f_\theta)^2 w_\theta - \sigma_{\xi,2}^2 w_\theta) \star K_{n,C_n}(Z)|^2]$$
$$\leq 3n^{-1}\{\mathbb{E}[(Y^4 + \sigma_{\xi,2}^4)|w_\theta \star K_{n,C_n}(Z)|^2] + 4\mathbb{E}[Y^2|(f_\theta w_\theta) \star K_{n,C_n}(Z)|^2]$$
$$+ \mathbb{E}[|(f_\theta^2 w_\theta) \star K_{n,C_n}(Z)|^2]\}$$
$$\leq \left(\frac{C}{n}\right)\{\mathbb{E}[f_{\theta^0}^4(X)(w_\theta \star K_{n,C_n})^2(Z)] + (\mathbb{E}(\xi^4) + \sigma_{\xi,2}^4)\mathbb{E}[(w_\theta \star K_{n,C_n})^2(Z)]$$
$$+ \mathbb{E}[f_{\theta^0}^2(X)((f_\theta w_\theta) \star K_{n,C_n}(Z))^2] + \mathbb{E}[((f_\theta^2 w_\theta) \star K_{n,C_n}(Z))^2]\},$$

that is, according to Lemma A.2

$$\mathrm{Var}[S_{n,1}(\theta)] \leq \left(\frac{C}{n}\right)\{\langle((f_{\theta^0}^4 g) \star p_\varepsilon), (w_\theta \star K_{n,C_n})^2\rangle + (\mathbb{E}(\xi^4) + \sigma_{\xi,2}^4)\langle(g \star p_\varepsilon), (w_\theta \star K_{n,C_n})^2\rangle$$
$$+ \langle((f_{\theta^0}^2 g) \star p_\varepsilon), ((f_\theta w_\theta) \star K_{n,C_n})^2\rangle + \sigma_{\xi,2}^2\langle(g \star p_\varepsilon), ((f_\theta w_\theta) \star K_{n,C_n})^2\rangle$$
$$+ \langle(g \star p_\varepsilon), ((f_\theta^2 w_\theta) \star K_{n,C_n})^2\rangle\}.$$

On one hand, according to (8.4), under the assumption ($A_3$),

$$\mathrm{Var}[S_{n,1}(\theta)] \leq \left(\frac{C}{n}\right)\{[\|(f_{\theta^0}^4 g) \star p_\varepsilon\|_\infty + (\mathbb{E}(\xi^4) + \sigma_{\xi,2}^4)\|g \star p_\varepsilon\|_\infty]\|w_\theta \star K_{n,C_n}\|_2^2$$
$$+ [\|(f_{\theta^0}^2 g) \star p_\varepsilon\|_\infty + \sigma_{\xi,2}^2\|g \star p_\varepsilon\|_\infty]\|(f_\theta w_\theta) \star K_{n,C_n}\|_2^2$$
$$+ \|g \star p_\varepsilon\|_\infty \|(f_\theta^2 w_\theta) \star K_{n,C_n}\|_2^2\},$$

and we get that

$$\mathrm{Var}[S_{n,1}(\theta)] \leq C(\sigma_{\xi,2}^2, \sigma_{\xi,4}, f_\theta, w_\theta, p_\varepsilon) n^{-1}\left[\left\|(w_\theta)^*\frac{K_{C_n}^*}{p_\varepsilon^*}\right\|_2^2 + \left\|(f_\theta w_\theta)^*\frac{K_{C_n}^*}{p_\varepsilon^*}\right\|_2^2 + \left\|(f_\theta^2 w_\theta)^*\frac{K_{C_n}^*}{p_\varepsilon^*}\right\|_2^2\right]. \quad (8.7)$$

On the other hand, again applying (8.4) under the assumption ($A_3$), we have

$$\mathrm{Var}[S_{n,1}(\theta)] \leq \frac{C}{n}\{[\|(f_{\theta^0}^4 g) \star p_\varepsilon\|_1 + (\mathbb{E}(\xi^4) + \sigma_{\xi,2}^4)\|g \star p_\varepsilon\|_1]\|w_\theta \star K_{n,C_n}\|_\infty^2$$
$$+ [\|(f_{\theta^0}^2 g) \star p_\varepsilon\|_1 + \sigma_{\xi,2}^2\|g \star p_\varepsilon\|_1]\|(f_\theta w_\theta) \star K_{n,C_n}\|_\infty^2$$
$$+ \|g \star p_\varepsilon\|_1\|(f_\theta^2 w_\theta) \star K_{n,C_n}\|_\infty^2\},$$

and hence we get that

$$\mathrm{Var}[S_{n,1}(\theta)] \leq C(\sigma_{\xi,2}^2, \sigma_{\xi,4}, f_\theta, w_\theta, p_\varepsilon) n^{-1}\left[\left\|(w_\theta)^*\frac{K_{C_n}^*}{p_\varepsilon^*}\right\|_1^2 + \left\|(f_\theta w_\theta)^*\frac{K_{C_n}^*}{p_\varepsilon^*}\right\|_1^2 + \left\|(f_\theta^2 w_\theta)^*\frac{K_{C_n}^*}{p_\varepsilon^*}\right\|_1^2\right]. \quad (8.8)$$

By combining (8.7) and (8.8) and by applying Lemma A.1 in the Appendix, we get that under assumption ($A_3$), $\mathrm{Var}[S_{n,1}(\theta)]$ is bounded by

$$n^{-1}C(\sigma_\xi^2, \sigma_{\xi,4}, f_\theta, w_\theta, p_\varepsilon) \times \min\left\{\left\|(w_\theta)^*\frac{K_{C_n}^*}{p_\varepsilon^*}\right\|_2^2 + \left\|(f_\theta w_\theta)^*\frac{K_{C_n}^*}{p_\varepsilon^*}\right\|_2^2 + \left\|(f_\theta^2 w_\theta)^*\frac{K_{C_n}^*}{p_\varepsilon^*}\right\|_2^2,\right.$$
$$\left.\left\|(w_\theta)^*\frac{K_{C_n}^*}{p_\varepsilon^*}\right\|_1^2 + \left\|(f_\theta w_\theta)^*\frac{K_{C_n}^*}{p_\varepsilon^*}\right\|_1^2 + \left\|(f_\theta^2 w_\theta)^*\frac{K_{C_n}^*}{p_\varepsilon^*}\right\|_1^2\right\}.$$

In other words,

$$\mathrm{Var}[S_{n,1}(\theta)] \leq C(\sigma_\xi^2, \sigma_{\xi,4}, f_\theta, w_\theta, p_\varepsilon) \max[1, C_n^{2\alpha-2a+(1-\rho)_+(1-\rho)_-} e^{-2bC_n^r + 2\beta C_n^\rho}]/n,$$

and under (3.7), then (8.2) is proved.

It remains now to check that there exists two sequences $\rho_k$ and $\epsilon_k$ tending to 0, such that (8.1) holds. First write that

$$|S_{n,1}(\theta) - S_{n,1}(\theta')| = \left|\frac{2}{n}\sum_{i=1}^n Y_i[(f_{\theta'} w_{\theta'} - f_\theta w_\theta)] \star K_{n,C_n}(Z_i) - \frac{1}{n}\sum_{i=1}^n [(f_{\theta'}^2 w_{\theta'} - f_\theta^2 w_\theta)] \star K_{n,C_n}(Z_i)\right|$$

$$\left| -\frac{1}{n} \sum_{i=1}^{n} (Y_i^2 - \sigma_{\xi,2}^2)[w_{\theta'} - w_\theta] \star K_{n,C_n}(Z_i) \right|$$

$$= \left| \frac{2}{n} \sum_{i=1}^{n} (\theta' - \theta)^\top \left( \frac{\partial [Y_i f_\theta w_\theta + f_\theta^2 w_\theta - (Y_i^2 - \sigma_{\xi,2}^2) w_\theta]}{\partial \theta} \bigg|_{\theta=\bar\theta} \right) \star K_{n,C_n}(Z_i) \right|,$$

where $\|\bar\theta\|_{\ell^2} \leq \|\theta - \theta'\|_{\ell^2}$. It follows that for $\|\theta - \theta'\|_{\ell^2} \leq \rho_k$,

$$|S_{n,1}(\theta) - S_{n,1}(\theta')| \leq \frac{2\rho_k}{n} \left\| \sum_{i=1}^{n} \left( \frac{\partial [Y_i f_\theta w_\theta + f_\theta^2 w_\theta + (Y_i^2 - \sigma_{\xi,2}^2) w_\theta]}{\partial \theta} \bigg|_{\theta=\bar\theta} \right) \star K_{n,C_n}(Z_i) \right\|_{\ell^2}.$$

Hence (8.1) holds since for $C_n$ satisfying (3.7), by using the same arguments as for the proof of (8.2), we have that for all $\theta \in \Theta$

$$\mathbb{E}\left[ \left\| n^{-1} \sum_{i=1}^{n} \left( \frac{\partial [Y_i f_\theta w_\theta + f_\theta^2 w_\theta + (Y_i^2 - \sigma_{\xi,2}^2) w_\theta]}{\partial \theta} \right) \star K_{n,C_n}(Z_i) \right\|_{\ell^2}^2 \right] = \mathrm{O}(1),$$

and hence for all $k \in \mathbb{N}$, $\mathbb{E}[\sup_{\|\theta'-\theta\|_{\ell^2} \leq \rho_k} |S_{n,1}(\theta) - S_{n,1}(\theta')|] = \mathrm{O}(\rho_k)$ as $n \to \infty$. □

*8.2. Proofs of (2) of Theorem 6.1 and Corollary 6.1*

If we denote by $S_{n,1}^{(1)}$ and $S_{n,1}^{(2)}$ the first and second derivatives of $S_{n,1}(\theta)$ with respect to $\theta$, by using classical Taylor expansion based on the smoothness properties of $\theta \mapsto w_\theta f_\theta$ and the consistency of $\widehat\theta_1$, we obtain that

$$0 = S_{n,1}^{(1)}(\widehat\theta_1) = S_{n,1}^{(1)}(\theta^0) + S_{n,1}^{(2)}(\theta^0)(\widehat\theta_1 - \theta^0) + R_{n,1}(\widehat\theta_1 - \theta^0),$$

with $R_{n,1}$ defined by

$$R_{n,1} = \int_0^1 [S_{n,1}^{(2)}(\theta^0 + s(\widehat\theta_1 - \theta^0)) - S_{n,1}^{(2)}(\theta^0)] \, \mathrm{d}s. \tag{8.9}$$

This implies that

$$\widehat\theta_1 - \theta^0 = -[S_{n,1}^{(2)}(\theta^0) + R_{n,1}]^{-1} S_{n,1}^{(1)}(\theta^0). \tag{8.10}$$

Consequently, we have to check the four following points.

(i) $\mathbb{E}[(S_{n,1}^{(1)}(\theta^0) - S_{\theta^0,g}^{(1)}(\theta^0))(S_{n,1}^{(1)}(\theta^0) - S_{\theta^0,g}^{(1)}(\theta^0))^\top] = \mathrm{O}[\varphi_n \varphi_n^\top]$.
(ii) $\mathbb{E}[\|S_{n,1}^{(2)}(\theta^0) - S_{\theta^0,g}^{(2)}(\theta^0)\|_{\ell^2}^2] = \mathrm{o}(1)$.
(iii) $R_{n,1}$ defined in (8.9) satisfies $\mathbb{E}(\|R_{n,1}\|_{\ell^2}^2) = \mathrm{o}(1)$ as $n \to \infty$.
(iv) $\mathbb{E}\|\widehat\theta_1 - \theta^0\|_{\ell^2}^2 \leq 4\mathbb{E}[(S_{n,1}^{(1)}(\theta^0))^\top [(S_{\theta^0,g}^{(2)}(\theta^0))^{-1}]^\top (S_{\theta^0,g}^{(2)}(\theta^0))^{-1} S_{n,1}^{(1)}(\theta^0)] + \mathrm{o}(\varphi_n^2)$.

The rate of convergence of $\widehat\theta_1$ is thus given by the order of $S_{n,1}^{(1)}(\theta^0) - S_{\theta^0,g}^{(1)}(\theta^0) = S_{n,1}^{(1)}(\theta^0)$.

**Proof of (i).** Write that

$$S_{n,1}^{(1)}(\theta) = \frac{\partial}{\partial \theta} \left( \frac{1}{n} \sum_{i=1}^{n} [(Y_i - f_\theta)^2 - \sigma_{\xi,2}^2] w_\theta \star K_{n,C_n}(Z_i) - \mathbb{E}[((Y - f_\theta(X))^2 - \sigma_{\xi,2}^2) w_\theta(X)] \right)$$

$$= \frac{1}{n} \sum_{i=1}^{n} \left( \frac{\partial [((Y_i - f_\theta)^2 - \sigma_{\xi,2}^2) w_\theta]}{\partial \theta} \right) \star K_{n,C_n}(Z_i) - \mathbb{E}\left[ \frac{\partial [((Y - f_\theta(X))^2 - \sigma_{\xi,2}^2) w_\theta(X)]}{\partial \theta} \right].$$

*Study of the bias.* According to Lemma A.2, $\mathbb{E}[S_{n,1}^{(1)}(\theta^0)]$ is equal to

$$-2\mathbb{E}\left[f_{\theta^0}(X)\left(\frac{\partial(f_\theta w_\theta)}{\partial \theta} \star K_{C_n}(X) - \frac{\partial(f_\theta w_\theta)}{\partial \theta}(X)\right)\right] + \mathbb{E}\left[\frac{\partial(f_\theta^2 w_\theta)}{\partial \theta} \star K_{C_n}(X) - \frac{\partial(f_\theta^2 w_\theta)}{\partial \theta}(X)\right]$$

$$+ \mathbb{E}\left[f_{\theta^0}^2(X)\left(\frac{\partial w_\theta}{\partial \theta} \star K_{C_n}(X) - \frac{\partial w_\theta}{\partial \theta}(X)\right)\right],$$

that is

$$\mathbb{E}[S_{n,1}^{(1)}(\theta^0)] = -2\left\langle (f_{\theta^0} g)^*, \left(\frac{\partial(f_\theta w_\theta)}{\partial \theta}\bigg|_{\theta=\theta^0}\right)^*(K_{C_n}^* - 1)\right\rangle$$

$$+ \left\langle g^*, \left(\frac{\partial(f_\theta^2 w_\theta)}{\partial \theta}\bigg|_{\theta=\theta^0}\right)^*(K_{C_n}^* - 1)\right\rangle$$

$$+ \left\langle (f_{\theta^0}^2 g)^*, \left(\frac{\partial w_\theta}{\partial \theta}\bigg|_{\theta=\theta^0}\right)^*(K_{C_n}^* - 1)\right\rangle.$$

On one hand, according to (8.4), under the assumption $(A_3)$, the bias is bounded in the following way

$$\left|\mathbb{E}\left[\frac{\partial S_{n,1}(\theta)}{\partial \theta_j}\right]_{\theta=\theta^0}\right| \leq \pi^{-1}\|(f_{\theta^0} g)^*\|_2 \left\|\left(\frac{\partial(f_\theta w_\theta)}{\partial \theta_j}\bigg|_{\theta=\theta^0}\right)^*(K_{C_n}^* - 1)\right\|_2$$

$$+ (2\pi)^{-1}\|g^*\|_2 \left\|\left(\frac{\partial(f_\theta^2 w_\theta)}{\partial \theta_j}\bigg|_{\theta=\theta^0}\right)^*(K_{C_n}^* - 1)\right\|_2$$

$$+ (2\pi)^{-1}\|(f_{\theta^0}^2 g)^*\|_2 \left\|\left(\frac{\partial w_\theta}{\partial \theta_j}\bigg|_{\theta=\theta^0}\right)^*(K_{C_n}^* - 1)\right\|_2$$

and consequently,

$$\left|\mathbb{E}\left[\left(\frac{\partial S_{n,1}(\theta)}{\partial \theta_j}\bigg|_{\theta=\theta^0}\right)\right]\right| \leq C(f_{\theta^0}, w_{\theta^0}, p_\varepsilon)\left[\left\|\left(\frac{\partial(w_\theta)}{\partial \theta_j}\bigg|_{\theta=\theta^0}\right)^*(K_{C_n}^* - 1)\right\|_2 \quad (8.11)\right.$$

$$\left. + \left\|\left(\frac{\partial(w_\theta f_\theta)}{\partial \theta_j}\bigg|_{\theta=\theta^0}\right)^*(K_{C_n}^* - 1)\right\|_2 + \left\|\left(\frac{\partial(w_\theta f_\theta^2)}{\partial \theta_j}\bigg|_{\theta=\theta^0}\right)^*(K_{C_n}^* - 1)\right\|_2\right].$$

And, on the other hand, the bias can also be bounded in the following way

$$\left|\mathbb{E}\left[\frac{\partial S_{n,1}(\theta)}{\partial \theta_j}\right]_{\theta=\theta^0}\right| \leq \pi^{-1}\|(f_{\theta^0} g)^*\|_\infty \left\|\left(\frac{\partial(f_\theta w_\theta)}{\partial \theta_j}\bigg|_{\theta=\theta^0}\right)^*(K_{C_n}^* - 1)\right\|_1$$

$$+ (2\pi)^{-1}\|g^*\|_\infty \left\|\left(\frac{\partial(f_\theta^2 w_\theta)}{\partial \theta_j}\bigg|_{\theta=\theta^0}\right)^*(K_{C_n}^* - 1)\right\|_1$$

$$+ (2\pi)^{-1}\|(f_{\theta^0}^2 g)^*\|_\infty \left\|\left(\frac{\partial w_\theta}{\partial \theta_j}\bigg|_{\theta=\theta^0}\right)^*(K_{C_n}^* - 1)\right\|_1$$

and consequently,

$$\left|\mathbb{E}\left[\left(\frac{\partial S_{n,1}(\theta)}{\partial \theta_j}\bigg|_{\theta=\theta^0}\right)\right]\right| \leq C(f_{\theta^0}, w_{\theta^0}, p_\varepsilon)\left[\left\|\left(\frac{\partial(w_\theta)}{\partial \theta_j}\bigg|_{\theta=\theta^0}\right)^*(K_{C_n}^* - 1)\right\|_1 \quad (8.12)\right.$$

$$\left. + \left\|\left(\frac{\partial(w_\theta f_\theta)}{\partial \theta_j}\bigg|_{\theta=\theta^0}\right)^*(K_{C_n}^* - 1)\right\|_1 + \left\|\left(\frac{\partial(w_\theta f_\theta^2)}{\partial \theta_j}\bigg|_{\theta=\theta^0}\right)^*(K_{C_n}^* - 1)\right\|_1\right].$$

By combining (8.11) and (8.12), we get that

$$\left|\mathbb{E}\left[\left(\frac{\partial S_{n,1}(\theta)}{\partial \theta_j}\bigg|_{\theta=\theta^0}\right)\right]\right| \leq C(f_{\theta^0}, w_{\theta^0}, p_\varepsilon) \min[B_{n,j}^{[1]}(\theta^0), B_{n,j}^{[2]}(\theta^0)],$$

with $B_{n,j}^{[q]}(\theta^0)$, $q=1,2$ defined in Theorem 6.1.

Consequently, by applying Lemma A.1, we obtain that

$$\mathbb{E}^2\left[\left(\frac{\partial S_{n,1}(\theta)}{\partial \theta_j}\bigg|_{\theta=\theta^0}\right)\right] \leq C_{b_2}^{(1)}(f_{\theta^0}, w_{\theta^0}, p_\varepsilon, b, r) C_n^{-2a+(1-r)+(1-r)_-} \mathrm{e}^{-2bC_n^r}.$$

*Study of the variance.* For the variance term, it is easy to see that

$$\mathrm{Var}\left(\frac{\partial S_{n,1}(\theta)}{\partial \theta_j}\bigg|_{\theta=\theta^0}\right) \leq \frac{C}{n}\mathbb{E}\left[\left(\frac{\partial[-2Y_i f_\theta w_\theta + f_\theta^2 w_\theta + (Y_i^2 - \sigma_{\xi,2}^2)w_\theta]}{\partial \theta_j}\bigg|_{\theta=\theta^0}\right) \star K_{n,C_n}(Z_i)\right]^2,$$

that is, according to Lemma A.2, $\mathrm{Var}(\partial S_{n,1}(\theta)/\partial \theta_j|_{\theta=\theta^0})$ equals

$$\frac{C+\mathrm{o}(1)}{n}\Bigg\{\left\langle((f_{\theta^0}^2 + \sigma_{\xi,2}^2)g) \star p_\varepsilon, \left(\left(\frac{\partial (f_\theta w_\theta)}{\partial \theta_j}\bigg|_{\theta=\theta^0}\right) \star K_{n,C_n}\right)^2\right\rangle$$
$$+ \left\langle((f_{\theta^0}^4 + \sigma_{\xi,4} + 4f_{\theta^0}^2 \sigma_{\xi,2}^2 - \sigma_{\xi,2}^4 + 4f_{\theta^0}\sigma_{\xi,3})g) \star p_\varepsilon, \left(\left(\frac{\partial w_\theta}{\partial \theta_j}\bigg|_{\theta=\theta^0}\right) \star K_{n,C_n}\right)^2\right\rangle$$
$$+ \left\langle g \star p_\varepsilon, \left(\left(\frac{\partial (f_\theta^2 w_\theta)}{\partial \theta_j}\bigg|_{\theta=\theta^0}\right) \star K_{n,C_n}\right)^2\right\rangle\Bigg\}.$$

It follows that according to (8.4), $\mathrm{Var}(\partial S_{n,1}(\theta)/\partial \theta_j|_{\theta=\theta^0})$ is less than

$$\frac{C}{n}\Bigg\{[\sigma_{\xi,2}^2\|g \star p_\varepsilon\|_\infty + \|(f_{\theta^0}^2 g) \star p_\varepsilon\|_\infty]\left\|\left(\frac{\partial (f_\theta w_\theta)}{\partial \theta_j}\bigg|_{\theta=\theta^0}\right) \star K_{n,C_n}\right\|_2^2$$
$$+ \|((f_{\theta^0}^4 + \sigma_{\xi,4} + 4f_{\theta^0}^2 \sigma_{\xi,2}^2 - \sigma_{\xi,2}^4 + 4f_{\theta^0}\sigma_{\xi,3})g) \star p_\varepsilon\|_\infty \left\|\left(\frac{\partial w_\theta}{\partial \theta_j}\bigg|_{\theta=\theta^0}\right) \star K_{n,C_n}\right\|_2^2$$
$$+ \|g \star p_\varepsilon\|_\infty \left\|\left(\frac{\partial (f_\theta^2 w_\theta)}{\partial \theta_j}\bigg|_{\theta=\theta^0}\right) \star K_{n,C_n}\right\|_2^2\Bigg\},$$

that is

$$\mathrm{Var}\left(\frac{\partial S_{n,1}(\theta)}{\partial \theta_j}\bigg|_{\theta=\theta^0}\right) \leq \frac{C(\sigma_{\xi,2}^2, f_{\theta^0}, f_{\theta^0,j}^{(1)}, w_{\theta^0}, p_\varepsilon)}{n}\Bigg[\left\|\left(\frac{\partial (w_\theta)}{\partial \theta_j}\bigg|_{\theta=\theta^0}\right)^* \frac{K_{C_n}^*}{p_\varepsilon^*}\right\|_2^2 \quad (8.13)$$
$$+ \left\|\left(\frac{\partial (f_\theta w_\theta)}{\partial \theta_j}\bigg|_{\theta=\theta^0}\right)^* \frac{K_{C_n}^*}{p_\varepsilon^*}\right\|_2^2 + \left\|\left(\frac{\partial (f_\theta^2 w_\theta)}{\partial \theta_j}\bigg|_{\theta=\theta^0}\right)^* \frac{K_{C_n}^*}{p_\varepsilon^*}\right\|_2^2\Bigg].$$

And once again, according to (8.4), another bound the variance term can be obtained to get that $\mathrm{Var}(\partial S_{n,1}(\theta)/\partial \theta_j|_{\theta=\theta^0})$ is bounded by

$$\frac{C}{n}\Bigg\{[\sigma_{\xi,2}^2\|g \star p_\varepsilon\|_1 + \|(f_{\theta^0}^2 g) \star p_\varepsilon\|_1]\left\|\left(\frac{\partial (f_\theta w_\theta)}{\partial \theta_j}\bigg|_{\theta=\theta^0}\right) \star K_{n,C_n}\right\|_\infty^2$$
$$+ \|g \star p_\varepsilon\|_1 \left\|\left(\frac{\partial (f_\theta^2 w_\theta)}{\partial \theta_j}\bigg|_{\theta=\theta^0}\right) \star K_{n,C_n}\right\|_\infty^2$$

$$+ \|((f_{\theta^0}^4 + \sigma_{\xi,4} + 4f_{\theta^0}^2\sigma_{\xi,2}^2 - \sigma_{\xi,2}^4 + 4f_{\theta^0}\sigma_{\xi,3})g) \star p_\varepsilon\|_1 \left\| \left(\frac{\partial w_\theta}{\partial \theta_j}\bigg|_{\theta=\theta^0}\right) \star K_{n,C_n} \right\|_\infty^2 \bigg\},$$

that is

$$\mathrm{Var}\left(\frac{\partial S_{n,1}(\theta)}{\partial \theta_j}\bigg|_{\theta=\theta^0}\right) \leq \frac{C(\sigma_{\xi,2}^2, f_{\theta^0}, f_{\theta^0,j}^{(1)}, w_{\theta^0}, p_\varepsilon)}{n}\left[\left\|\left(\frac{\partial(w_\theta)}{\partial \theta_j}\bigg|_{\theta=\theta^0}\right)^* \frac{K_{C_n}^*}{p_\varepsilon^*}\right\|_1^2 \right. \quad (8.14)$$
$$+ \left\|\left(\frac{\partial(f_\theta w_\theta)}{\partial \theta_j}\bigg|_{\theta=\theta^0}\right)^* \frac{K_{C_n}^*}{p_\varepsilon^*}\right\|_1^2 + \left.\left\|\left(\frac{\partial(f_\theta^2 w_\theta)}{\partial \theta_j}\bigg|_{\theta=\theta^0}\right)^* \frac{K_{C_n}^*}{p_\varepsilon^*}\right\|_1^2\right].$$

By combining (8.13) and (8.14), we get that

$$\mathrm{Var}\left(\frac{\partial S_{n,1}(\theta)}{\partial \theta_j}\bigg|_{\theta=\theta^0}\right) \leq \frac{C(\sigma_{\xi,2}^2, f_{\theta^0}, f_{\theta^0,j}^{(1)}, w_{\theta^0}, p_\varepsilon)}{n} \min\{V_{n,j}^{[1]}(\theta^0), V_{n,j}^{[2]}(\theta^0)\},$$

with $V_{n,j}^{[q]}$, $q = 1, 2$ defined in Theorem 6.1. By applying Lemma A.1, we get that

$$\mathrm{Var}\left(\frac{\partial S_{n,1}(\theta)}{\partial \theta_j}\bigg|_{\theta=\theta^0}\right) \leq \frac{C_{v_2}(\sigma_{\xi,2}^2, f_{\theta^0}f_{\theta^0,j}^{(1)}, w_{\theta^0}, p_\varepsilon)}{n} \max[1, C_n^{2\alpha-2a+(1-\rho)+(1-\rho)_-}\exp\{-2bC_n^r + 2\beta C_n^\rho\}].$$

The rate of convergence of $\widehat{\theta}_1$ denoting by $\varphi_n \varphi_n^\top$ corresponds to the best choice for the sequence $C_n^*$, minimizing the sum of the variance $\mathrm{Var}[S_{n,1}^{(1)}(\theta^0)]$ and the square of the bias $(\mathbb{E}[S_{n,1}^{(1)}(\theta^0)] - S^{(1)}(\theta^0))(\mathbb{E}[S_{n,1}^{(1)}(\theta^0)] - S^{(1)}(\theta^0))^\top = (\mathbb{E}[S_{n,1}^{(1)}(\theta^0)])(\mathbb{E}[S_{n,1}^{(1)}(\theta^0)])^\top$.

*Polynomial noise* (see (N$_2$) with $\beta = \rho = 0$).

- If for $j = 1, \ldots, m$, $\partial(f_\theta w_\theta)/\partial\theta_j|_{\theta=\theta^0}$, $\partial w_\theta/\partial\theta_j|_{\theta=\theta^0}$ and $\partial(f_\theta^2 w_\theta)/\partial\theta_j|_{\theta=\theta^0}$ satisfy (R$_1$) with $r = 0$, then

$$\left|\mathbb{E}\left[\frac{\partial S_{n,1}(\theta)}{\partial \theta_j}\bigg|_{\theta=\theta^0}\right]\right|^2 \leq C_{b_2}(f_{\theta^0}, w_{\theta^0}, b, \sigma_{\xi,2}^2)C_n^{-2a+1},$$

and

$$\mathrm{Var}\left(\frac{\partial S_{n,1}(\theta)}{\partial \theta_j}\bigg|_{\theta=\theta^0}\right) \leq C_{v_2}(\mathbb{E}(f_{\theta^0}^2(X)), \sigma_{\xi,2}^2, f_{\theta^0}^{(1)}, f_{\theta^0}, w_{\theta^0}, p_\varepsilon)\max[1, C_n^{2\alpha-2a+1}]/n.$$

It follows that if $r = 0$ and $a < \alpha + 1/2$ then

$$C_n^* = n^{1/(2\alpha)} \quad \text{and} \quad \varphi_n^2 = \mathrm{O}(n^{(1-2a)/(2\alpha)}). \qquad (8.15)$$

If $r = 0$ and $a \geq \alpha + 1/2$, then

$$C_n^* = n^{1/(2a-1)} \quad \text{and} \quad \varphi_n^2 = \mathrm{O}(n^{-1}).$$

- If for $j = 1, \ldots, m$, $\partial(f_\theta w_\theta)/\partial\theta_j|_{\theta=\theta^0}$, $\partial w_\theta/\partial\theta_j|_{\theta=\theta^0}$ and $\partial(f_\theta^2 w_\theta)/\partial\theta_j|_{\theta=\theta^0}$ satisfy (R$_1$) with $r > 0$, then

$$\left|\mathbb{E}\left[\frac{\partial S_{n,1}(\theta)}{\partial \theta_j}\bigg|_{\theta=\theta^0}\right]\right|^2 \leq C_{b_2}(f_{\theta^0}, \sigma_{\xi,2}^2, f_{\theta^0}^{(1)}, w_{\theta^0})C_n^{-2a+(1-r)+(1-r)_-}\exp\{-2bC_n^r\}$$

and

$$\mathrm{Var}\left(\frac{\partial S_{n,1}(\theta)}{\partial \theta_j}\bigg|_{\theta=\theta^0}\right) \leq \frac{C_{v_2}(f_{\theta^0}, \sigma_{\xi,2}^2, f_{\theta^0}^{(1)}, w_{\theta^0}, p_\varepsilon)}{n}.$$

It follows that

$$C_n^* = \left[\frac{\log n}{2b} + \frac{-2a + (1-r) + (1-r)_-}{2br}\log\left(\frac{\log n}{2b}\right)\right]^{1/r} \quad \text{and} \quad \varphi_n^2 = \mathrm{O}(n^{-1}). \tag{8.16}$$

*Exponential noise* (see (N$_2$) with $\rho > 0$).

• If for $j = 1, \ldots, m$, $\partial(f_\theta w_\theta)/\partial\theta_j|_{\theta=\theta^0}$, $\partial w_\theta/\partial\theta_j|_{\theta=\theta^0}$ and $\partial(f_\theta^2 w_\theta)/\partial\theta_j|_{\theta=\theta^0}$ satisfy (R$_1$) with $r = 0$, then

$$\left|\mathbb{E}\left[\frac{\partial S_{n,1}(\theta)}{\partial \theta_j}\bigg|_{\theta=\theta^0}\right]\right|^2 \leq C_{b_2}(f_{\theta^0}, \sigma_{\xi,2}^2, f_{\theta^0}^{(1)}, w_{\theta^0}) C_n^{-2a+1},$$

and

$$\mathrm{Var}\left(\frac{\partial S_{n,1}(\theta)}{\partial \theta_j}\bigg|_{\theta=\theta^0}\right) \leq C_{v_2}(f_{\theta^0}, \sigma_{\xi,2}^2, f_{\theta^0}^{(1)}, w_{\theta^0}, p_\varepsilon) C_n^{2\alpha - 2a + (1-\rho) + (1-\rho)_-} \frac{\exp\{2\beta C_n^\rho\}}{n}.$$

It follows that

$$C_n^* = \left[\frac{\log n}{2\beta} - \frac{2\alpha + (1-\rho)_-}{2\rho\beta}\log\left(\frac{\log n}{2\beta}\right)\right]^{1/\rho} \quad \text{and} \quad \varphi_n^2 = \mathrm{O}\left[\left(\frac{\log n}{2\beta}\right)^{(1-2a)/\rho}\right]. \tag{8.17}$$

• If for $j = 1, \ldots, m$, $\partial(f_\theta w_\theta)/\partial\theta_j|_{\theta=\theta^0}$, $\partial w_\theta/\partial\theta_j|_{\theta=\theta^0}$ and $\partial(f_\theta^2 w_\theta)/\partial\theta_j|_{\theta=\theta^0}$ satisfy (R$_1$) with $r > 0$ and $\{r > \rho\}$ or $\{r = \rho$ and $b > \beta\}$ or $\{r = \rho$, $b = \beta$ and $a \geq \alpha + 1/2\}$ then

$$\left|\mathbb{E}\left[\frac{\partial S_{n,1}(\theta)}{\partial \theta_j}\bigg|_{\theta=\theta^0}\right]\right|^2 \leq C_{b_2}(\mathbb{E}(f_{\theta^0}^2(X)), \sigma_{\xi,2}^2, L(f_{\theta^0}^{(1)} w_{\theta^0}), L(f_{\theta^0} f_{\theta^0}^{(1)} w_{\theta^0})) C_n^{-2a+(1-r)+(1-r)_-} \exp\{-2bC_n^r\},$$

and

$$\mathrm{Var}\left(\frac{\partial S_{n,1}(\theta)}{\partial \theta_j}\bigg|_{\theta=\theta^0}\right) \leq C_{v_2}(\mathbb{E}(f_{\theta^0}^2(X)), \sigma_{\xi,2}^2, L(f_{\theta^0}^{(1)} w_{\theta^0}), L(f_{\theta^0} f_{\theta^0}^{(1)} w_{\theta^0}))/n.$$

It follows that

$$C_n^* = \left[\frac{\log n}{2b} + \frac{-2a + (1-r) + (1-r)_-}{2br}\log\left(\frac{\log n}{2b}\right)\right]^{1/r} \quad \text{and} \quad \varphi_n^2 = \mathrm{O}(n^{-1}). \tag{8.18}$$

• If for $j = 1, \ldots, m$, $\partial(f_\theta w_\theta)/\partial\theta_j|_{\theta=\theta^0}$, $\partial w_\theta/\partial\theta_j|_{\theta=\theta^0}$ and $\partial(f_\theta^2 w_\theta)/\partial\theta_j|_{\theta=\theta^0}$ satisfy (R$_1$) with $r > 0$ and $r \leq \rho$ or $r = \rho$ and $b < \beta$ then

$$\left|\mathbb{E}\left[\frac{\partial S_{n,1}(\theta)}{\partial \theta_j}\bigg|_{\theta=\theta^0}\right]\right|^2 \leq C_b(\mathbb{E}(f_{\theta^0}^2(X)), \sigma_{\xi,2}^2, L(f_{\theta^0}^{(1)} w_{\theta^0}), L(f_{\theta^0} f_{\theta^0}^{(1)} w_{\theta^0})) C_n^{-2a+(1-r)+(1-r)_-} \exp\{-2bC_n^r\}$$

and

$$\mathrm{Var}\left(\frac{\partial S_{n,1}(\theta)}{\partial \theta_j}\bigg|_{\theta=\theta^0}\right) \leq C_{v_2}(\mathbb{E}(f_{\theta^0}^2(X)), \sigma_{\xi,2}^2, L(f_{\theta^0}^{(1)} w_{\theta^0}), L(f_{\theta^0} f_{\theta^0}^{(1)} w_{\theta^0}))$$

$$\times C_n^{2\alpha - 2a + (1-\rho) + (1-\rho)_-} \exp\{-2bC_n^r + 2\beta C_n^\rho\}/n.$$

It follows that the bias is by a logarithmic factor larger than the variance and

$$C_n^* = \left[\frac{\log n}{2\beta} - \frac{2\alpha + (1-\rho)_- - (1-r)_-}{2\rho\beta}\log\left(\frac{\log n}{2\beta}\right)\right]^{1/\rho}$$

and

$$\varphi_n^2 = \mathrm{O}\left[(\log n)^{(-2a+(1-r)+(1-r)_-)/\rho} \exp\left\{-2b\left(\frac{\log n}{2\beta}\right)^{r/\rho}\right\}\right].$$

If $r = \rho$, $b = \beta$ and $a < \alpha + 1/2$ then

$$\mathrm{Var}\left(\frac{\partial S_{n,1}(\theta)}{\partial \theta_j}\bigg|_{\theta=\theta^0}\right) \leq C_{v_2}(\mathbb{E}(f_{\theta^0}^2(X)), \sigma_{\xi,2}^2, L(f_{\theta^0}^{(1)}w_{\theta^0}), L(f_{\theta^0}f_{\theta^0}^{(1)}w_{\theta^0}))\frac{C_n^{2\alpha-2a+1}}{n}.$$

The resulting rate is thus $\varphi_n^2 = \mathrm{O}[(\log n)^{(2\alpha-2a+1)/r}/n]$. $\square$

**Proof of (ii).** By using that

$$(S_{n,1}^{(2)}(\theta))_{j,k} = -\frac{2}{n}\sum_{i=1}^{n}\left(-2Y_i\frac{\partial^2(f_\theta w_\theta)}{\partial\theta_j\partial\theta_k} + \frac{\partial^2(f_\theta^2 w_\theta)}{\partial\theta_j\partial\theta_k} + (Y_i^2 - \sigma_{\xi,2}^2)\frac{\partial^2 w_\theta}{\partial\theta_j\partial\theta_k}\right)\bigg|_{\theta=\theta^0} \star K_{n,C_n}(Z_i), \qquad (8.19)$$

we write $S_{n,1}^{(2)}(\theta^0) - S_{\theta^0,g}^{(2)}(\theta^0) = A_0 + A_1 + A_2 + A_3$ with

$$A_0 = \mathbb{E}(S_{n,1}^{(2)}(\theta^0)) - S_{\theta^0,g}^{(2)}(\theta^0),$$

$$A_1 = -\frac{2}{n}\sum_{i=1}^{n}Y_i\left(\frac{\partial^2(f_\theta w_\theta)}{\partial\theta_j\partial\theta_k}\bigg|_{\theta=\theta^0}\right)\star K_{n,C_n}(Z_i) - \mathbb{E}\left[Y_i\left(\frac{\partial^2(f_\theta w_\theta)}{\partial\theta_j\partial\theta_k}(X)\bigg|_{\theta=\theta^0}\right)\right],$$

$$A_2 = \frac{1}{n}\sum_{i=1}^{n}\left(\frac{\partial^2(f_\theta^2 w_\theta)}{\partial\theta_j\partial\theta_k}\bigg|_{\theta=\theta^0}\right)\star K_{n,C_n}(Z_i) - \mathbb{E}\left[\left(\frac{\partial^2(f_\theta^2 w_\theta)}{\partial\theta_j\partial\theta_k}(X)\bigg|_{\theta=\theta^0}\right)\right],$$

$$A_3 = \frac{1}{n}\sum_{i=1}^{n}(Y_i^2 - \sigma_{\xi,2}^2)\left(\frac{\partial^2 w_\theta}{\partial\theta_j\partial\theta_k}\bigg|_{\theta=\theta^0}\right)\star K_{n,C_n}(Z_i) - \mathbb{E}\left[(Y_i^2 - \sigma_{\xi,2}^2)\left(\frac{\partial^2 w_\theta}{\partial\theta_j\partial\theta_k}(X)\bigg|_{\theta=\theta^0}\right)\right].$$

As soon as $w_\theta f_\theta$ and $w_\theta f_\theta^2$ and their derivatives up to order 2, satisfy $(R_1)$, then for $C_n$ satisfying that (3.7), we get that $A_0 = \mathrm{o}(1)$ and $A_j = \mathrm{o}_p(1)$ as $n \to \infty$ for $j = 1, \ldots, 3$, and (ii) is proved. $\square$

**Proof of (iii).** Again using (8.19), the smoothness properties of the derivatives of $w_\theta f_\theta$ and $w_\theta f_\theta^2$ up to order 3 and the consistency of $\widehat{\theta}_1$, we get that $\|R_{n,1}\|_{\ell^2} = \mathrm{o}_p(1)$ as $n \to \infty$. $\square$

**Proof of (iv).** Let us introduce the random event $E_n = \bigcap_{j,k} E_{n,j,k}$, where

$$E_{n,j,k} = \left\{\omega \text{ such that } \left|\frac{\partial^2 S_{\theta^0,g}(\theta)}{\partial\theta_j\partial\theta_k}\bigg|_{\theta=\theta^0} - \frac{\partial^2 S_{n,1}(\theta,\omega)}{\partial\theta_j\partial\theta_k}\bigg|_{\theta=\theta^0} + (R_{n,1})_{j,k}(\omega)\right| \leq \frac{1}{2}\frac{\partial^2 S_{\theta^0,g}(\theta)}{\partial\theta_j\partial\theta_k}\bigg|_{\theta=\theta^0}\right\}.$$

According to (8.10) and (8.9), we have

$$\mathbb{E}\|\widehat{\theta}_1 - \theta^0\|_{\ell^2}^2 = \mathbb{E}[\|\widehat{\theta}_1 - \theta^0\|_{\ell^2}^2 \mathbb{1}_{E_n}] + \mathbb{E}[\|\widehat{\theta}_1 - \theta^0\|_{\ell^2}^2 \mathbb{1}_{E_n^c}]$$

$$\leq \mathbb{E}[\|\widehat{\theta}_1 - \theta^0\|_{\ell^2}^2 \mathbb{1}_{E_n}] + 2\sup_{\theta\in\Theta}\|\theta\|_{\ell^2}^2 \mathbb{P}_{\theta^0,g}(E_n^c)$$

$$\leq \mathbb{E}[(S_{n,1}^{(1)}(\theta^0))^\top [(S_{n,1}^{(2)}(\theta^0) + R_{n,1})^{-1}]^\top (S_{n,1}^{(2)}(\theta^0) + R_{n,1})^{-1} S_{n,1}^{(1)}(\theta^0)\mathbb{1}_{E_n}]$$

$$+ 2\sup_{\theta\in\Theta}\|\theta\|_{\ell^2}^2 \mathbb{P}_{\theta^0,g}(E_n^c)$$

$$\leq C 4^{2m^2} \sup_{j,k}\left|\frac{\partial^2 S_{\theta^0,g}(\theta)}{\partial\theta_j\partial\theta_k}\bigg|_{\theta=\theta^0}\right|^{-2} \mathbb{E}[(S_{n,1}^{(1)}(\theta^0))^\top S_{n,1}^{(1)}(\theta^0)] + 2\sup_{\theta\in\Theta}\|\theta\|_{\ell^2}^2 \mathbb{P}_{\theta^0,g}(E_n^c)$$

$$\leq C 4^{2m^2} \sup_{j,k} \left| \frac{\partial^2 S_{\theta^0,g}(\theta)}{\partial \theta_j \partial \theta_k} \bigg|_{\theta=\theta^0} \right|^{-2} \varphi_n^2 + 2 \sup_{\theta \in \Theta} \|\theta\|_{\ell^2}^2 \mathbb{P}_{\theta^0,g}(E_n^c).$$

It remains thus to show that $\mathbb{P}_{\theta^0,g}(E_n^c) = \mathrm{o}(\varphi_n^2)$ with

$$\sup_{j,k} \mathbb{E}\left[\left(\frac{\partial^2 (S_{n,1}(\theta) - S_{\theta^0,g}(\theta))}{\partial \theta_j \partial \theta_k}\bigg|_{\theta=\theta^0}\right)^2\right] \leq \varphi_n^2.$$

By Markov's inequality, for a $p > 2$,

$$\mathbb{P}_{\theta^0,g}(E_n^c) \leq \sum_{j=1}^m \sum_{k=1}^m \mathbb{P}_{\theta^0,g}(E_{n,j,k}^c),$$

with

$$\mathbb{P}_{\theta^0,g}(E_{n,j,k}^c) \leq \mathbb{E}\left[\left|\left(\frac{\partial^2 (S_{\theta^0,g}(\theta) - S_{n,1}(\theta))}{\partial \theta_j \partial \theta_k}\bigg|_{\theta=\theta^0}\right) + (R_{n,1})_{j,k}\right|^p\right]$$

$$\leq 2^{p-1} \left|\left(\frac{\partial^2 S_{\theta^0,g}(\theta)}{\partial \theta_j \partial \theta_k}\bigg|_{\theta=\theta^0}\right) - \mathbb{E}\left[\left(\frac{\partial^2 S_{n,1}(\theta)}{\partial \theta_j \partial \theta_k}\bigg|_{\theta=\theta^0}\right)\right]\right|^p$$

$$+ 2^{p-1}\mathbb{E}\left[\left|\mathbb{E}\left(\frac{\partial^2 S_{n,1}(\theta)}{\partial \theta_j \partial \theta_k}\bigg|_{\theta=\theta^0}\right) - \left(\frac{\partial^2 S_{n,1}(\theta)}{\partial \theta_j \partial \theta_k}\bigg|_{\theta=\theta^0}\right) + (R_{n,1})_{j,k}\right|^p\right].$$

In other words,

$$\mathbb{P}_{\theta^0,g}(E_{n,j,k}^c) \leq 2^{p-1} \left|\left(\frac{\partial^2 S_{\theta^0,g}(\theta)}{\partial \theta_j \partial \theta_k}\bigg|_{\theta=\theta^0}\right) - \mathbb{E}\left[\left(\frac{\partial^2 S_{n,1}(\theta)}{\partial \theta_j \partial \theta_k}\bigg|_{\theta=\theta^0}\right)\right]\right|^p$$

$$+ 2^{2p-2}\left\{\mathbb{E}\left[\left|\mathbb{E}\left(\frac{\partial^2 S_{n,1}(\theta)}{\partial \theta_j \partial \theta_k}\bigg|_{\theta=\theta^0}\right) - \left(\frac{\partial^2 S_{n,1}(\theta)}{\partial \theta_j \partial \theta_k}\bigg|_{\theta=\theta^0}\right)\right|^p\right] + 2^{2p-2}\mathbb{E}|(R_{n,1})_{j,k}|^p\right\}.$$

Now we apply the Rosenthal's inequality to the sum of centered variables

$$\left(\frac{\partial^2 S_{n,1}(\theta)}{\partial \theta_j \partial \theta_k}\bigg|_{\theta=\theta^0}\right) - \mathbb{E}\left[\left(\frac{\partial^2 S_{n,1}(\theta)}{\partial \theta_j \partial \theta_k}\bigg|_{\theta=\theta^0}\right)\right] = n^{-1}\sum_{i=1}^n W_{n,i,j,k},$$

where $W_{n,i,j,k}$ equals

$$\left(\frac{\partial^2 [-2Y_i f_\theta w_\theta + f_\theta^2 w_\theta + (Y_i^2 - \sigma_{\xi,2}^2)w_\theta]}{\partial \theta_j \partial \theta_k}\bigg|_{\theta=\theta^0}\right) \star K_{n,C_n}(Z_i)$$

$$- 2\mathbb{E}\left[\left(\frac{\partial^2 [-2Y_i f_\theta w_\theta + f_\theta^2 w_\theta + (Y_i^2 - \sigma_{\xi,2}^2)w_\theta]}{\partial \theta_j \partial \theta_k}\bigg|_{\theta=\theta^0}\right) \star K_{n,C_n}(Z_i)\right].$$

**Lemma 8.1 (Rosenthal's inequality [27]).** *For $U_1, \ldots, U_n$, be $n$ independent centered random variables, there exists a constant $C(p)$ such that for $p \geq 1$,*

$$\mathbb{E}\left[\left|\sum_{i=1}^n U_i\right|^p\right] \leq C(p)\left[\sum_{i=1}^n \mathbb{E}[|U_i|^p] + \left(\sum_{i=1}^n \mathbb{E}[U_i^2]\right)^{p/2}\right]. \tag{8.20}$$

Consequently,

$$\mathbb{E}\left[\left|\left(\frac{\partial^2 S_{n,1}(\theta)}{\partial \theta_j \, \partial \theta_k}\bigg|_{\theta=\theta^0}\right) - \mathbb{E}\left(\left(\frac{\partial^2 S_{n,1}(\theta)}{\partial \theta_j \, \partial \theta_k}\bigg|_{\theta=\theta^0}\right)\right)\right|^p\right] \le C(p)[n^{1-p}\mathbb{E}|W_{n,1,j,k}|^p + n^{-p/2}\mathbb{E}^{p/2}|W_{n,1,j,k}|^2].$$

Take $p = 4$ to get that

$$\mathbb{E}\left[\left|\left(\frac{\partial^2 S_{n,1}(\theta)}{\partial \theta_j \, \partial \theta_k}\bigg|_{\theta=\theta^0}\right) - \mathbb{E}\left(\left(\frac{\partial^2 S_{n,1}(\theta)}{\partial \theta_j \, \partial \theta_k}\bigg|_{\theta=\theta^0}\right)\right)\right|^4\right] \le C(4)[n^{-3}\mathbb{E}|W_{n,1}|^4 + n^{-2}\mathbb{E}^2|W_{n,1}|^2].$$

Therefore, under the conditions ensuring that $(\mathbb{E}[S^{(2)}_{\theta^0,g}(\theta^0) - S^{(2)}_{n,1}(\theta^0)]^2)_{j,k} = \mathrm{o}(1)$, we have

$$(\mathbb{E}[S^{(2)}_{n,1}(\theta^0) - S^{(2)}_{\theta^0,g}(\theta^0)]^4)_{j,k} = \mathrm{O}(\varphi_n^4) = \mathrm{o}(\varphi_n^2).$$

Now, by using the definition of $R_{n,1}$ combined with (8.19) and the smoothness properties of the derivatives of $(w_\theta f_\theta)$ and $(w_\theta f_\theta^2)$ up to order 3, we get that $\mathbb{E}((R_{n,1})^4_{j,k}) = \mathrm{o}(\|\widehat{\theta}_1 - \theta^0\|^4_{\ell^2})$, and we conclude that

$$\mathbb{E}\|\widehat{\theta}_1 - \theta^0\|^2_{\ell^2} \le 4\mathbb{E}[(S^{(1)}_{n,1}(\theta^0))^\top [(S^{(2)}_{\theta^0,g}(\theta^0))^{-1}]^\top (S^{(2)}_{\theta^0,g}(\theta^0))^{-1} S^{(1)}_{n,1}(\theta^0)]$$
$$+ \mathrm{o}(\varphi_n^2) + \mathrm{o}(\mathbb{E}[\|\widehat{\theta}_1 - \theta^0\|^4_{\ell^2}]). \qquad \Box$$

### 8.3. Proof of Theorem 6.2

The proof of Theorem 6.2 follows from Theorem 6.1 from the conditions $(C_8)$–$(C_{10})$ and from the central limit theorem with the Lindeberg condition (see, for instance, [4]). The main point of the proof lies in the proof of

$$\sqrt{n} S^{(1)}_{n,1}(\theta^0) \xrightarrow[n \to \infty]{\mathcal{L}} \mathcal{N}(0, \Sigma_{0,1}), \tag{8.21}$$

where $\Sigma_{0,1}$ is defined in Theorem 6.2. We start by writing that

$$\sqrt{n} S^{(1)}_{n,1}(\theta^0) = \sqrt{n}[S^{(1)}_{n,1}(\theta^0) - S^{(1)}_{\theta^0,g}(\theta^0)]$$
$$= \sqrt{n}[S^{(1)}_{n,1}(\theta^0) - \mathbb{E}(S^{(1)}_{n,1}(\theta^0))] + \sqrt{n}[\mathbb{E}(S^{(1)}_{n,1}(\theta^0)) - S^{(1)}_{\theta^0,g}(\theta^0)]$$
$$= \sum_{i=1}^n \frac{V_{n,i}}{\sqrt{n}} + \sqrt{n}[\mathbb{E}(S^{(1)}_{n,1}(\theta^0)) - S^{(1)}_{\theta^0,g}(\theta^0)].$$

where

$$V_{n,i} = \frac{\partial[((Y_i - f_\theta)^2 - \sigma^2_{\xi,2}) w_\theta]}{\partial \theta}\bigg|_{\theta=\theta^0} \star K_{n,C_n}(Z_i) - \mathbb{E}\left(\frac{\partial[((Y_i - f_\theta)^2 - \sigma^2_{\xi,2}) w_\theta]}{\partial \theta}\bigg|_{\theta=\theta^0} \star K_{n,C_n}(Z_i)\right).$$

The proof of (8.21) follows from the three points

(A-1) $n^{-1} \sum_{i=1}^n V_{n,i} V_{n,i}^\top \xrightarrow[n \to \infty]{\mathbb{P}} \widetilde{\Sigma}_{0,1}$;

(A-2) For all $\epsilon > 0$, $n^{-1} \sum_{i=1}^n V_{n,i} V_{n,i}^\top \mathbb{1}_{\|V_{n,i}\|_{\ell^2} > \epsilon \sqrt{n}} \xrightarrow[n \to \infty]{\mathbb{P}} 0$;

(A-3) $\sqrt{n}[\mathbb{E}(\widetilde{S}^{(1)}_{n,1}(\theta^0)) - S^{(1)}_{\theta^0,g}(\theta^0)] \xrightarrow[n \to \infty]{} 0$.

**Proof of (A-1).** The conditions $(C_8)$–$(C_{10})$ ensure that $V_{n,j} = \mathrm{O}(1)$ for $j = 1, \ldots, d$. Hence, (A-1) is checked. $\qquad \Box$

**Proof of (A-2).** It suffices to prove that for all $\epsilon > 0$ and for all $j = 1, \ldots, d$,

$$\mathbb{E}\left[n^{-1}\sum_{i=1}^{n}(V_{n,i,j})^2 \mathbb{1}_{|V_{n,i,j}|>\epsilon\sqrt{n}}\right] = o(1),$$

where $V_{n,i,j}$ is the $j$th coordinate of the vector $V_{n,i}$. Now (A-2) is checked by writing that under (A$_5$) and (C$_8$)–(C$_{10}$),

$$n^{-1}\sum_{i=1}^{n}\mathbb{E}[(V_{n,i,j})^2 \mathbb{1}_{|V_{n,i,j}|>\epsilon\sqrt{n}}] \leq \frac{1}{n\sqrt{n}\epsilon}\sum_{i=1}^{n}\mathbb{E}|V_{n,i,j}|^3$$

$$\leq \frac{1}{\sqrt{n}\epsilon} \times O(1). \qquad \square$$

**Proof of (A-3).** Conditions (C$_8$)–(C$_{10}$) ensure that for all $C_n$ the variance term $V_{n,j} = O(1)$ for $j = 1, \ldots, d$. Hence, it is possible to choose $C_n$ such that (A-3) is checked. $\qquad \square$

*8.4. Proof of Theorem 6.3*

As for the proof of Theorem 6.1, the main point of the proof consists in showing that for any $\theta \in \Theta$, $\mathbb{E}[(S_{n,2}(\theta) - S_{\theta^0,g}(\theta))^2] = o(1)$, with $S_{\theta^0,g}(\theta)$ admitting a unique minimum in $\theta = \theta^0$. The second part of the proof consists in studying $\omega_1(n, \rho)$ defined as $\omega_1(n, \rho) = \sup\{|S_{n,2}(\theta) - S_{n,2}(\theta')|: \|\theta - \theta'\|_2 \leq \rho\}$. By using the regularity assumptions on the regression function $f$, we state that there exists two sequences $\rho_k$ and $\epsilon_k$ tending to 0, such that for all $k \in \mathbb{N}$,

$$\lim_{n\to\infty}\mathbb{P}[\omega_1(n,\rho_k) > \epsilon_k] = 0 \quad \text{and that} \quad \mathbb{E}[(\omega_1(n,\rho_k))^2] = O(\rho_k^2). \tag{8.22}$$

Under the conditions (C$_{11}$)-(C$_{13}$) and by applying the law of large numbers, we get that for any $\theta \in \Theta$, $\mathbb{E}[(S_{n,2}(\theta) - S_{\theta^0,g}(\theta))^2] = o(1)$ as $n \to \infty$ and $\mathbb{E}[(S_{n,2}(\theta) - S_{n,2}(\theta'))^2] = o(1)$ and consequently, $\widehat{\theta}_2$ is a consistent estimator of $\theta^0$. By using classical Taylor expansion based on the smoothness properties of the regression function, with respect to $\theta$ and the consistency of $\widehat{\theta}_2$, we get that $0 = S_{n,2}^{(1)}(\widehat{\theta}_1) = S_{n,2}^{(1)}(\theta^0) + S_{n,2}^{(2)}(\theta^0)(\widehat{\theta}_2 - \theta^0) + R_{n,2}(\widehat{\theta}_2 - \theta^0)$, with $R_{n,2}$ defined by

$$R_{n,2} = \int_0^1 [S_{n,2}^{(2)}(\theta^0 + s(\widehat{\theta}_2 - \theta^0))\widehat{\theta} - S_{n,2}^{(2)}(\theta^0)]\,ds. \tag{8.23}$$

This implies that $\widehat{\theta}_2 - \theta^0 = -[S_{n,2}^{(2)}(\theta^0) + R_{n,2}]^{-1}S_{n,2}^{(1)}(\theta^0)$. The $\sqrt{n}$-consistency and the asymptotic normality follow by applying the central limit theorem with the Lindeberg condition to $S_{n,2}^{(1)}(\theta^0)$ to get that $\sqrt{n}S_{n,2}^{(1)}(\theta^0) \xrightarrow[n\to\infty]{\mathcal{L}} \mathcal{N}(0, \Sigma_{0,2})$.

**Appendix: Technical lemmas**

**Lemma A.1.** *Let $\varphi$ a function such that $\varphi$ belongs to $\mathbb{L}_2(\mathbb{R})$ satisfying (R$_1$). Then*

$$\int_{|u|\geq C_n}|\varphi^*(u)|\,du \leq \frac{L(\varphi)}{R(a,b,r)}C_n^{-a+1-r}\exp\{-bC_n^r\}. \tag{A.1}$$

*Furthermore, if $p_\varepsilon$ satisfies (N$_2$), then*

$$\int_{|u|\leq C_n}\frac{|\varphi^*(u)|}{|p_\varepsilon^*(u)|}\,du \leq \frac{L(\varphi)}{R(\alpha,\beta,\rho,a,b,r)\underline{C}(p_\varepsilon)}\max[1, C_n^{\alpha-a+1-\rho}\exp\{-bC_n^r + \beta C_n^\rho\}].$$

**Lemma A.2.** *Let $\varphi$ such that $\mathbb{E}(|\varphi(Y,X)|)$ is finite and let $\Phi$ such that $\mathbb{E}(|\Phi(U)|)$ is finite. Then*

$$\mathbb{E}[\varphi(Y,X)\Phi \star K_{n,C_n}(U)] = \mathbb{E}[\varphi(Y,X)\Phi \star K_{C_n}(X)] = \langle \varphi(y,\cdot)f_{Y,X}(y,\cdot), \Phi \star K_{C_n}\rangle$$

*and*

$$\mathbb{E}[\varphi(Y,X)\Phi \star K_{n,C_n}(U)]^2 = \int \langle (\varphi(x,\cdot)f_{Y,X}(x,\cdot)) \star f_\varepsilon, (\Phi \star K_{n,C_n})^2\rangle \, dx.$$

**Proof.** If we denote by $f_{Y,Z,X}$ the joint distribution of $(Y,Z,X)$, then

$$f_{Y,Z,X}(y,z,x) = f_{Y,X}(y,x)f_\varepsilon(z-x) \quad \text{and} \quad f_{Y,Z}(y,z) = f_{Y,X}(y,\cdot) \star f_\varepsilon(z).$$

By construction $K_{n,C_n} \star p_\varepsilon = K_{C_n}$, and hence Parseval's formula

$$\begin{aligned}
\mathbb{E}[\varphi(Y,X)\Phi \star K_{n,C_n}(Z)] &= \iint \varphi(y,x)f_{Y,X}(y,x) \int \Phi \star K_{n,C_n}(z)f_\varepsilon(z-x)\,dz\,dy\,dx \\
&= (2\pi)^{-1} \iint \varphi(y,x)f_{Y,X}(y,x) \int \Phi^*(u)\frac{K^*_{C_n}(u)}{f^*_\varepsilon(u)}f^*_\varepsilon(u)e^{-iux}\,du\,dy\,dx \\
&= \iint \varphi(y,x)f_{Y,X}(y,x) \int \Phi(z)K_{C_n}(x-z)\,dz\,dy\,dx.
\end{aligned}$$

The second equality follows by writing that

$$\begin{aligned}
\mathbb{E}[\varphi(Y,X)\Phi \star K_{n,C_n}(Z)]^2 &= \iiint \varphi(y,x)(\Phi \star K_{n,C_n}(z))^2 f_{Y,X}(y,x)f_\varepsilon(z-x)\,dy\,dz\,dx \\
&= \iiint \varphi(y,x)(\Phi \star K_{n,C_n}(z))^2 f_{Y,X}(y,x)f_\varepsilon(z-x)\,dy\,dz\,dx \\
&= \int \langle (\varphi(y,\cdot)f_{Y,X}(y,\cdot)) \star f_\varepsilon, (\Phi \star K_{n,C_n})^2\rangle\,dy.
\end{aligned}\qquad\square$$

## References

[1] S. Baran. A consistent estimator in general functional errors-in-variables models. *Metrika* **51** (2000) 117–132 (electronic). MR1790927
[2] P. J. Bickel, A. J. C. Klaassen, Y. Ritov and J. A. Wellner. *Efficient and Adaptative Estimation for Semiparametric Model*. Johns Hopkins Univ. Press, Baltimore, MD, 1993. MR1245941
[3] Bickel, P. J. and A. J. C. Ritov. Efficient estimation in the errors-in-variables model. *Ann. Statist.* **15** (1987) 513–540. MR0888423
[4] Billingsley, P. *Probability and Measure*, 3rd edition. Wiley. New York, 1995. MR1324786
[5] R. J. Carroll, D. Ruppert and L. A. Stefanski. *Measurement Error in Nonlinear Models*. Chapman and Hall, London, 1995. MR1630517
[6] L. K. Chan and T. K. Mak. On the polynomial functionnal relationship. *J. Roy. Statist. Soc. Ser. B* **47** (1985) 510–518. MR0844482
[7] C. H. Cheng and J. W. Van Ness. On estimating linear relationships when both variables are subject to errors. *J. Roy. Statist. Soc. Ser. B* **56** (1994) 167–183. MR1257805
[8] F. Comte and M.-L. Taupin. Semiparametric estimation in the (auto)-regressive $\beta$-mixing model with errors-in-variables. *Math. Methods Statist.* **10** (2001) 121–160. MR1851745
[9] I. Fazekas, S. Baran, A. Kukush, and J. Lauridsen. Asymptotic properties in space and time of an estimator in nonlinear functional errors-in-variables models. *Random Oper. Stochastic Equations* **7** (1999) 389–412. MR1709899
[10] I. Fazekas and A. G. Kukush. Asymptotic properties of estimators in nonlinear functional errors-in-variables with dependent error terms. *J. Math. Sci. (New York)* **92** (1998) 3890–3895. MR1666219
[11] M. V. Fedoryuk. *Asimptotika: integraly i ryady*. "Nauka", Moscow, 1987. MR0950167
[12] W. A. Fuller. *Measurement Error Models*. Wiley, New York, 1987. MR0898653
[13] L. J. Gleser. Improvements of the naive approach to estimation in nonlinear errors-in-variables regression models. *Contemp. Math.* **112** (1990) 99–114. MR1087101

[14] J. A. Hausman, W. K. Newey, I. Ichimura and J. L. Powell. Identification and estimation of polynomial errors-in-variables models. *J. Econometrics* **50** (1991) 273–295. MR1147115

[15] J. A. Hausman, W. K. Newey and J. L. Powell. Nonlinear errors in variables estimation of some engel curves. *J. Econometrics* **65** (1995) 205–233. MR1324193

[16] H. Hong and E. Tamer. A simple estimator for nonlinear error in variable models. *J. Econometrics* **117** (2003) 1–19. MR2002282

[17] C. Hsiao. Consistent estimation for some nonlinear errors-in-variables models. *J. Econometrics* **41** (1989) 159–185. MR1007729

[18] C. Hsiao, L. Wang and Q. Wang. Estimation of nonlinear errors-in-variables models: an approximate solution. *Statist. Papers* **38** (1997) 1–25. MR1474937

[19] C. Hsiao and Q. K. Wang. Estimation of structural nonlinear errors-in-variables models by simulated least-squares method. *Internat. Econom. Rev.* **41** (2000) 523–542. MR1760462

[20] J. Kiefer and J. Wolfowitz. Consistency of the maximum likelihood estimator in the presence of infinitely many nuisance parameters. *Ann. Math. Statist.* **27** (1956) 887–906. MR0086464

[21] A. Kukush and H. Schneeweiss. Comparing different estimators in a nonlinear measurement error model. I. *Math. Methods Statist.* **14** (2005) 53–79. MR2158071

[22] A. Kukush and H. Schneeweiss. Comparing different estimators in a nonlinear measurement error model. II. *Math. Methods Statist.* **14** (2005) 203–223. MR2160395

[23] O. V. Lepski and B. Y. Levit. Adaptive minimax estimation of infinitely differentiable functions. *Math. Methods Statist.* **7** (1998) 123–156. MR1643256

[24] T. Li. Estimation of nonlinear errors-in-variables models: a simulated minimum distance estimator. *Statist. Probab. Lett.* **47** (2000) 243–248. MR1747484

[25] T. Li. Robust and consistent estimation of nonlinear errors-in-variables models. *J. Econometrics* **110** (2002) 1–26. MR1920960

[26] S. A. Murphy and A. W. Van der Vaart. Likelihood inference in the errors-in-variables model. *J. Multivariate Anal.* **59** (1996) 81–108. MR1424904

[27] V. V. Petrov. *Limit Theorems of Probability Theory*. Oxford Science Publications, New York, 1995. MR1353441

[28] O. Reiersøl. Identifiability of a linear relation between variables which are subject to error. *Econometrica.* **18** (1950) 375–389. MR0038054

[29] M.-L. Taupin. Semi-parametric estimation in the nonlinear structural errors-in-variables model. *Ann. Statist.* **29** (2001) 66–93. MR1833959

[30] A. van der Vaart. Semiparametric statistics. *Lectures on Probability Theory and Statistics (Saint-Flour, 1999)* 331–457. *Lecture Notes in Math.* **1781**. Berlin, Springer, 2002. MR1915446

[31] A. W. van der Vaart. Estimating a real parameter in a class of semiparametric models. *Ann. Statist.* **16** (1988) 1450–1474. MR0964933

[32] A. W. van der Vaart. Efficient estimation in semi-parametric mixture models. *Ann. Statist.* **24** (1996) 862–878. MR1394993

[33] K. M. Wolter and W. A. Fuller. Estimation of nonlinear errors-in variables models. *Ann. Statist.* **10** (1982) 539–548. MR0653528

[34] K. M. Wolter and W. A. Fuller. Estimation of the quadratic errors-in-variables model. *Biometrika* **69** (1982) 175–182.